\documentclass[12pt]{amsart}
\usepackage{amssymb,amsmath,graphicx}
\makeatletter
\def\@cite#1#2{{\m@th\upshape\bfseries%
[{#1\if@tempswa{\m@th\upshape\mdseries, #2}\fi}]}}
\makeatother
\theoremstyle{plain}
\newtheorem{thm}{Theorem}[section]
\newtheorem{lem}[thm]{Lemma}
\newtheorem{cor}[thm]{Corollary}
\newtheorem{prop}[thm]{Proposition}
\theoremstyle{definition}
\newtheorem{rem}[thm]{Remark}
\newtheorem{defn}[thm]{Definition}
\newtheorem{note}[thm]{Note}
\newtheorem{eg}[thm]{Example}

\newtheorem{const}[thm]{Construction}
\theoremstyle{remark}
\newtheorem*{acknowledgements}{Acknowledgements}

\newlength{\qedskip}

\newcommand{\BH}{{\B(\H)}}
\newcommand{\ca}{\mathit{C}^*}

\newcommand{\ol}{\overline}
\newenvironment{sbmatrix}{\left[\begin{smallmatrix}}{\end{smallmatrix}\right]}

\newcommand{\td}{\widetilde}

\newcommand{\bbC}{{\mathbb{C}}}

\newcommand{\bbF}{{\mathbb{F}}}

\newcommand{\bbR}{{\mathbb{R}}}
\newcommand{\bbT}{{\mathbb{T}}}
\newcommand{\bbZ}{{\mathbb{Z}}}
 \newcommand{\A}{{\mathcal{A}}}
 \newcommand{\B}{{\mathcal{B}}}

 \newcommand{\E}{{\mathcal{E}}}
 \newcommand{\F}{{\mathcal{F}}}
 
\renewcommand{\H}{{\mathcal{H}}}

 \newcommand{\K}{{\mathcal{K}}}  

 \newcommand{\M}{{\mathcal{M}}}
 \newcommand{\N}{{\mathcal{N}}}
\renewcommand{\O}{{\mathcal{O}}}

\renewcommand{\S}{{\mathcal{S}}}
 \newcommand{\T}{{\mathcal{T}}}
 
 \newcommand{\V}{{\mathcal{V}}}

\renewcommand{\phi}{\varphi}
\newcommand{\upchi}{{\raise.35ex\hbox{$\chi$}}}
\newcommand{\fA}{{\mathfrak{A}}}

\newcommand{\qand}{\quad\text{and}\quad}

\newcommand{\qfor}{\quad\text{for}\quad}

\newcommand{\ran}{\operatorname{Ran}}

\newcommand{\fockspaceN}{\operatorname{\F}\nolimits_N}
\newcommand{\fockspacetwoN}{\operatorname{\F}\nolimits_{2N}}

\newcommand{\uhfNI}{\operatorname{UHF}\nolimits_{N^\infty}}

\newcommand{\sumiN}{\sum_{i=1}^N}

\newcommand{\bofh}{\B(\H)}

\newbox\ipbox
\newcommand{\ip}[2]{\left\langle #1\mathrel{\mathchoice
{\setbox\ipbox=\hbox{$\displaystyle \left\langle\mathstrut #1#2\right\rangle$}
\vrule height\ht\ipbox width0.25pt depth\dp\ipbox}
{\setbox\ipbox=\hbox{$\textstyle \left\langle\mathstrut #1#2\right\rangle$}
\vrule height\ht\ipbox width0.25pt depth\dp\ipbox}
{\setbox\ipbox=\hbox{$\scriptstyle \left\langle\mathstrut #1#2\right\rangle$}
\vrule height\ht\ipbox width0.25pt depth\dp\ipbox}
{\setbox\ipbox=\hbox{$\scriptscriptstyle \left\langle\mathstrut #1#2\right\rangle$}
\vrule height\ht\ipbox width0.25pt depth\dp\ipbox}
} #2\right\rangle}
\newcommand{\tS}{\smash{\tilde{S}}\vphantom{S}}
\newcommand{\tA}{\smash{\tilde{A}}\vphantom{A}}
\newcommand{\tT}{\smash{\tilde{T}}\vphantom{T}}
\newcommand{\tm}{\smash{\tilde{m}}\vphantom{m}}

\begin{document}

\title[Wavelet representations and Fock-space constructions]%
{Wavelet representations \\ and  Fock space
on positive matrices\rlap{${}^{1}$}}
\author[P.E.T. Jorgensen]{P.E.T.~Jorgensen${}^{2}$}
\author[D.W.Kribs]{D.W.~Kribs\rlap{${}^{3}$}}
\thanks{\llap{${}^{1\;}$}The research in this paper was completed in the
Fall of 2000, but for various reasons the paper was late coming out in its
final form.
\newline\indent\leavevmode\llap{${}^{2\;}$}The first author was
partially supported by the National
Science Foundation.
\newline\indent\leavevmode\llap{${}^{3\;}$}The second author was partially
supported
by a Canadian NSERC
Post-doctoral Fellowship.%
}

\subjclass%
{42C40, 42A16, 43A65, 42A65}

\begin{abstract}
We show that every biorthogonal wavelet determines
a representation by operators on Hilbert space satisfying
simple identities, which captures
the established relationship between orthogonal wavelets
and Cuntz-algebra representations in that special case.
Each of these representations is shown to have tractable finite-%
dimensional co-invariant doubly-cyclic subspaces.
Further, motivated by these representations, we
introduce a general Fock-space Hilbert space construction which
yields creation operators containing the Cuntz--%
Toeplitz isometries as a special case.
\end{abstract}
\maketitle


In this paper, we wish to establish a connection
between biorthogonal wavelets on the one hand \cite{Dau92},
and representation theory for operators on Hilbert space
on the other \cite{BJendII,DKS}. This is accomplished by
showing that each of these wavelets yields a
collection of operators acting on Hilbert space
which satisfy simple identities, and which contain the
Cuntz relations \cite{Cun} as a special case. In fact,
this new relationship collapses to the now
well-known connection between orthogonal 
wavelets and representations of the Cuntz $\ca$-algebra
in that special case \cite{BJwave}.
Our second goal is to develop a framework for
studying this new class of representations.
Toward this end, we introduce a general
Fock space Hilbert space construction which reduces to unrestricted
Fock space in the familiar cases.
Indeed, the natural creation operators we get can be thought of
as
an
analogue of the Cuntz--Toeplitz creation operators to this more general
setting. We regard this construction and the 
creation operators  determined by it
as interesting objects of study in their own right. Finally,
our hope is that this paper will lead to
further study of the relationships and objects introduced
here.

\section{Introduction}\label{S:int}

In recent years we have seen several operator-theoretic approaches
to wavelet theory, e.g.,
\cite{BaLa99,BJwave,HLPS99,Jorgen_min,Kri_cp,PaRi01}.
Typically they involve representing wavelets of a particular type by
operators on infinite-dimensional Hilbert space. They have had success
because often the operators satisfy simple identities and hence lend
themselves to investigation.  From an operator theory and operator algebra
point of view, these
approaches often yield interesting new classes of examples to work with,
and can open up new areas of study
\cite{BrJo01,DKS,DP1,Jor99a,Jorgen_min}. From the wavelet perspective,
operator theory can be used to study the
fundamental analysis/synthesis problem for wavelets, i.e., the
transformation from functions in the Hilbert space $L^2(\bbR)$ to wavelet
coefficients in $\ell^2(\bbZ)$. For instance, the paper \cite{Jorgen_min}
includes an application of work from \cite{DKS} on {\it free semigroup
operator algebras} to obtain a lucid characterization of when the data
going into an {\it orthogonal wavelet} is minimal. 
Pioneering early papers which suggest use of operator theory
and representations of groups and algebras in wavelet analysis include
\cite{BaLa99},
\cite{BaMe99},
\cite{BGRW99},
\cite{CDF92},
\cite{Dau95},
\cite{HLPS99},
\cite{Kla99},
\cite{Law90},
\cite{PaRi01}, and
\cite{Wic93}.
In this paper, we introduce an operator-theoretic approach for the study
of {\it biorthogonal wavelets}. We will discuss the particulars of
such wavelets and this approach in the next section. 

We now discuss some of the basics of the wavelet cum operator approach,
and we will continue this discussion in the next section. We use standard
wavelet nomenclature from such texts as \cite{BrJo01,Dau92}. Wavelet
theory is centred around the action of the integers $\bbZ$ on the Hilbert
space $L^2(\bbR)$ by the translations $f\mapsto f\left( x-k\right) $, 
$k\in\mathbb{Z}$, and by a fixed {\it scaling operator}
\begin{equation}
U\colon f\longmapsto\frac{1}{\sqrt{N}}f\left( \frac{x}{N}\right)  \qfor
f\in L^2(\bbR). 
\label{Scaling}
\end{equation}
The existence of a {\it resolution subspace} $\V$ in $L^2(\bbR)$ which is
invariant for translation by $\bbZ$ and also for the scaling operator $U$,
is equivalent to the existence of a cyclic subspace under translation with
cyclic
vector $\varphi \in L^2(\bbR)$ which in an important special case 
can be shown to satisfy 
\begin{equation}
\frac{1}{\sqrt{N}} \, \varphi \big(\frac{x}{N} \big) = 
\sum_{k\in \bbZ} a_k \, \varphi (x-k),
\end{equation}
for some scalars $a_k$. When such a $\varphi$ exists, it is called a {\it
scaling function} and generates in an algorithmic fashion a set of wavelet
basis functions for $L^2(\bbR)$. 
There is also an explicit correspondence between the scaling function and
its so called {\it filter functions} $\{ m_i: 0 \leq i \leq N-1\}$, which
in turn define the operators $\{S_i \} $ we have been studying. For the
biorthogonal
wavelets, there will be two functions $\varphi$, $\td{\varphi}$, as
well as two sets of related filter functions $\{ m_i, \td{m}_j\}$ and sets
of operators $\{ S_i, \td{S}_j \}$ which we are interested in. One of the
attractive properties of this whole set up is that the various
correspondences are explicit; there are formulae which allow us to go back
and forth between each of the settings. We expand on these correspondences
below, and in the next section.

Just as passing to a resolution subspace $L^{2}\left( \mathbb{R}\right)
\rightarrow\mathcal{V}$ is a reduction of the initial
analysis/synthesis problem,
we will aim for a setup which is effective for computations,
and so a Hilbert space isometry
$\mathcal{V}\rightarrow\ell^{2}\left( \mathbb{Z}\right) $ is
desirable, and a further reduction to a much smaller subspace,
which converts the original problem into one of
manipulating sequences, i.e., vectors in $\ell^{2}\left( \mathbb{Z}\right) $.
But we
have $\ell^{2}\left( \mathbb{Z}\right) \simeq L^{2}\left( \mathbb{T}\right)$,
$\mathbb{T}=\mathbb{R}\diagup 2\pi\mathbb{Z}$, by virtue of
the Fourier transform, and it will be convenient to
couch the operator theory in terms of the \emph{function space}
$L^{2}\left( \mathbb{T}\right)$. Setting $N=2$ for simplicity, it turns out
that it is possible to find functions $m_{0}$, $m_{1}$ on $\mathbb{T}$, and
functions $\varphi$, $\psi$ in $L^{2}\left( \mathbb{R}\right)$
with $\varphi\in\mathcal{V}$ and $U\psi\in\mathcal{V}$,
such that the
quadrature wavelet problem takes the following form: Let
$S_{i}h\left( z\right) =m_{i}\left( z\right) h\left( z^{2}\right) $,
$h\in L^{2}\left( \mathbb{T}\right)$, $i=0,1$, and let $\hat{S}_{i}$ be the
corresponding operators on $\ell^{2}\left( \mathbb{Z}\right) $ with
adjoints $\hat{S}_{i}^{\ast}$. Introduce for $c=\{c_k \}_{k\in\bbZ} \in
\ell^2(\bbZ)$, 
\begin{equation}
c\ast f\left( x\right) =\sum_{k\in\mathbb{Z}}c_{k}f\left( x-k\right) ,
\qquad x\in\mathbb{R}.
\label{Composition}
\end{equation}
Then under suitable conditions on $m_{0}$ and $m_{1}$, it is
possible to get the following representation of the scaling
operator (called the resolution/detail
decomposition):
\[
\mathcal{V}\ni c\ast\varphi =U\left\lbrack 
\left( \hat{S}_{0}^{\ast}c\right) \ast\varphi\right\rbrack 
+U\left\lbrack \left( \hat{S}_{1}^{\ast}c\right) \ast\psi\right\rbrack ,
\]
where $U$ is the scaling operator (\ref{Scaling}) for $N=2$. From
this we can then deduce an algorithmic approach to
the analysis/synthesis problem of wavelets, i.e.,
\[
L^{2}\left( \mathbb{R}\right) \ni f=\sum_{j,k\in\mathbb{Z}}c_{j,k}\psi_{j,k},
\qquad f\longleftrightarrow\left( c_{j,k}\right) ,
\]
where $\psi_{j,k}\left( x\right) 
=2^{\frac{j}{2}}\psi\left( 2^{j}x-k\right) $ is a wavelet
basis for $L^{2}\left( \mathbb{R}\right) $. 

Hence the wavelet problem has
been translated into one for a different operator
system, not in $L^{2}\left( \mathbb{R}\right) $ but rather in the sequence
space $\ell^{2}\left( \mathbb{Z}\right) $. The operators $F_i =
\hat{S}_i^*$ are known in signal processing as subband filters. When they
are further assumed to satisfy the identities $(i)$ of
Theorem~\ref{ProS:wavelets.2}, they are called quadrature mirror filters.
The identities $(i)$ are also called the Cuntz relations. Realizing that
the Cuntz relations are in fact the quadrature subband-filter identities,
therefore makes the connection to the so called `pyramid algorithm' of
wavelet theory, see \cite{Dau92}; in other words, to the problem of
discretizing signals using wavelets.

The next reduction is then to
try to determine the latter problem from an equivalent one which
involves only a finite-dimensional subspace $\mathcal{K}$ in
$\ell^{2}\left( \mathbb{Z}\right) $, or equivalently in
$L^{2}\left( \mathbb{T}\right) $. In other words, the goal is to discern
the actions of the $S_i$-operators, hence the structure of the resolution
and
the wavelet itself, simply by examining their actions on a tractable
finite-dimensional subspace. This has been accomplished for orthogonal
wavelets when the scaling function $\varphi$ has compact support.
The space $\mathcal{K}$ is called the \emph{anchor subspace} and it has
been studied in separate and independent earlier papers by the coauthors
\cite{BrJo01,DKS,Jorgen_min,Kri_cp}. We will show here that this can be
done effectively for biorthogonal wavelets as well.



The next section begins with a discussion of the
general method used to represent orthogonal wavelets as operators on
Hilbert space satisfying the Cuntz relations. 
In particular, we recall the equivalence between
(i)~orthogonal wavelets of scale $N$,
(ii)~certain representations of the Cuntz $\ca$-algebra $\O_N$, and
(iii)~\emph{matricial loops}, i.e., the
group of all bounded measurable
functions from the torus $\mathbb{T}$ into
the unitary group $\mathrm{U}_{N}\left( \mathbb{C}\right) $.
We then prove a corresponding result for biorthogonal wavelets,  including
an analogous matrix perspective
which involves \emph{invertible loops}, i.e., the larger
non-compact group of all bounded measurable
functions from $\mathbb{T}$ into the
general linear group ($=\mathrm{GL}_{N}\left( \mathbb{C}\right) $.)
The link between a wavelet and its loop is provided by the 
filter functions. 

As for the orthogonal wavelet representations
\cite{DKS,Jorgen_min,Kri_cp}, every biorthogonal wavelet 
representation is shown to have tractable finite-dimensional
co-invariant cyclic subspaces.
The finite dimensionality of the subspace requires that the wavelet
be of compact support.
This is the content
of the third section. In particular, the representation
can be recovered, in a spatial sense, from these finite-dimensional anchor
subspaces. We remark on how this can be regarded as a weak dilation
theorem for the operators determined by the biorthogonal wavelet. 
We also find a striking relationship between the
invertible-loop matrix entries and the operators $\{ S_i, \td{S}_j\}$ 
which determine a biorthogonal representation. This gives us motivation
for the Fock space approach. 

In the final two sections, we introduce a general
Fock space Hilbert space construction which is motivated by the
biorthogonal wavelet
representations. Every completely positive map from complex
matrices to the bounded operators on a Hilbert space
(equivalently, every positive matrix with operator entries)
determines a `twisted' Fock-space structure,
which in turn yields natural {\it creation operators}.
The Cuntz--Toeplitz isometries acting on unrestricted
Fock space are recovered in
two different ways; through the
orthogonal class, and then from how the orthogonal class
sits inside the biorthogonal class. We describe the
actions of the creation operators which the spatial
construction yields, and characterize when they are bounded
operators in terms of the completely positive map.
Our initial goal with this construction was to find
analogues of the Cuntz--Toeplitz isometries for the
biorthogonal wavelet representations. In any event,
we regard these new creation operators as
interesting objects of study in their own right, and
our hope is that this paper will lead to further
study of them.

\section{Wavelet representations and operator identities}\label{S:wavelets}

We shall consider a family
of representations associated
with wavelets, and relate them
to representations of operator
identities, which have the
Cuntz relations as a special
case. In the simplest
case, the associated operators on $L^2\left(\mathbb{T}\right)$ in the
representations have the form
\begin{equation}
S\colon f\longmapsto m\left(z\right)f\left(z^{N}\right),
\qquad f\in L^{2}\left(\mathbb{T}\right),\label{eqS:wavelets.1}
\end{equation}
where $\mathbb{T}=\mathbb{R}\diagup 2\pi\mathbb{Z}$ is the
usual torus, and $L^{2}\left(\mathbb{T}\right)$ is the
Hilbert space defined from the Haar
measure $\mu$ on $\mathbb{T}$. The number $N
\in\left\{2,3,4,\dots\right\}$ will be fixed,
and the function $m\in L^{\infty}\left(\mathbb{T}\right)$
determines the operator.

Using the isomorphism
$\ell^2\left(\mathbb{Z}\right)\cong L^2\left(\mathbb{T}\right)$
of Fourier series, we note that
the operator (\ref{eqS:wavelets.1}) may also be
realized as acting on sequences $\mathbf{x}$ as follows:
\begin{equation}
\left(\mathop{S}\mathbf{x}\right)_i
=\sum_{j\in\mathbb{Z}}c_{i-Nj}x_j,
\text{\qquad for }i\in\mathbb{Z}
\text{ and }\mathbf{x}=\left(x_j\right)_{j\in\mathbb{Z}},
\label{eqS:wavelets.2}
\end{equation}
where $\{c_j\}_{j\in\bbZ}$ forms the Fourier expansion of $m(z)$. 
The corresponding $\infty\times\infty$ matrix
has the following form (the case $N=2$):
\[
\begin{array}
[c]{c}%
\setlength{\unitlength}{12pt}
\begin{picture}(16,16)(-8,-8)
\put(-7.6,0.6){\includegraphics
[bb=0.357178 0.309449 14.6428 12.6904,width=15.2\unitlength,height=13bp]
{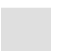}}
\put(-1.4,-5.4){\includegraphics
[bb=0.357178 0.309449 14.6428 12.6904,width=26bp,height=10.7\unitlength]
{shade.eps}}
\put(0,0){\makebox(0,0){$\displaystyle\begin{array}{ccccccccc}
 &  &  &  & \vdots &  &  &  &  \\
 & c_{3} & c_{1}
 & \multicolumn{1}{|c}{c_{-1}}
 & c_{-3} & c_{-5} & c_{-7} & c_{-9} &  \\
 & c_{4} & c_{2}
 & \multicolumn{1}{|c}{c_{0}}
 & c_{-2} & c_{-4} & c_{-6} & c_{-8} &  \\\cline{4-4}
 & c_{5} & c_{3}
 & c_{1} & \multicolumn{1}{|c}{c_{-1}}
 & c_{-3} & c_{-5} & c_{-7} &  \\
\cdots
 & c_{6} & c_{4}
 & c_{2} & \multicolumn{1}{|c}{c_{0}}
 & c_{-2} & c_{-4} & c_{-6} & 
\cdots \\\cline{5-5}
 & c_{7} & c_{5}
 & c_{3} & c_{1} & \multicolumn{1}{|c}{c_{-1}}
 & c_{-3} & c_{-5} &  \\
 & c_{8} & c_{6}
 & c_{4} & c_{2} & \multicolumn{1}{|c}{c_{0}}
 & c_{-2} & c_{-4} &  \\\cline{6-6}
 & c_{9} & c_{7}
 & c_{5} & c_{3} & c_{1} & \multicolumn{1}{|c}{c_{-1}}
 & c_{-3} &  \\
 & c_{10} & c_{8}
 & c_{6} & c_{4} & c_{2} & \multicolumn{1}{|c}{c_{0}}
 & c_{-2} &  \\\cline{7-7}
 & c_{11} & c_{9}
 & c_{7} & c_{5} & c_{3} & c_{1} & \multicolumn{1}{|c}{c_{-1}}
 &  \\
 &  &  &  & \vdots &  &  &  & 
\end{array}$}}\end{picture}%
\end{array}
\]
It is down-slanted with slope $2$, and it
is called a \emph{slanted Toeplitz matrix}. Its
spectral properties are in some ways analogous to
those of Toeplitz matrices, and in other ways
completely different; see \cite{BrJo01} and \cite{CDM91}.
In this form, it is known as the
\emph{subdivision operator}, and it is
used in numerical analysis, see
e.g., \cite{BrJo01}, \cite{CDM91}, \cite{Kla99}, \cite{Wic93}.
The translation from (\ref{eqS:wavelets.1}) to (\ref{eqS:wavelets.2}) may
be carried out by the usual Fourier series representation,
\begin{equation}
m\left(z\right)=\sum_{k\in\mathbb{Z}}c_kz^k,
\qquad c_{k}=\int_{\mathbb{T}}z^{-k}m\left( z\right) \,d\mu\left( z\right) ,
\label{FourierSeries}
\end{equation}
and
\[
f\left(z\right)=\sum_{k\in\mathbb{Z}}x_kz^k,\qquad z\in\mathbb{T},
\; x_{k}=\int_{\mathbb{T}}z^{-k}f\left( z\right) \,d\mu\left( z\right) .
\]
The numbers $\left(c_k\right)$ are called the
\textit{masking coefficients}
for
the subdivision.
The adjoint of (\ref{eqS:wavelets.1}) is known as the
transfer operator, alias the
Perron--Frobenius--Ruelle operator \cite{BrJo01}, and its action is
readily described. The proof of the following lemma is routine. 

\begin{lem}
\label{LemActionAdjoint}
If $S$ is an operator on $L^{2}\left(\mathbb{T}\right)$ with
$\left(Sf\right)\left(z\right)=m\left(z\right)f\left(z^{N}\right)$,
then $S^{*}$ is given by
\[
\left(S^{*}f\right)\left(z\right)=\frac{1}{N}\sum_{w^{N}=z}
\overline{m\left(w\right)}\,f\left(w\right).
\]
\end{lem}


Let us expand on the connection between the operators $S$ on $L^2
(\mathbb{T})$ and the underlying function theory on $L^2(\bbR)$. 

\begin{prop}
Let $\varphi \in L^2(\bbR)$ satisfy $(2)$. Define the operator $W_\varphi
: L^2(\bbT) \rightarrow L^2(\bbR)$ by 
\[
\big( W_\varphi f \big) (x) = \sum_{k\in \bbZ} x_k \varphi (x-k), \,\,\,\,
x_k = \hat{f}(k) = \int_{\bbT} z^{-k}f(z) \, d\mu (z), \,\,\, k\in\bbZ. 
\]
Then
\begin{equation}
\frac{1}{\sqrt{N}} \big( W_\varphi f \big) \big( \frac{x}{N}\big) =
\big(W_\varphi S f \big) (x), \qfor x\in\bbR, 
\end{equation}
where $S$ is the operator defined by $(4)$ and $(6)$. 
\end{prop}

\begin{proof}
Let $\varphi$, $W_\varphi$ and $S$ be as stated. Then for $f(z) =
\sum_{k\in\bbZ} x_k z^k$, we conclude
\begin{eqnarray*}
\frac{1}{\sqrt{N}} \big(W_\varphi f \big) \big( \frac{x}{N} \big) &=&
\frac{1}{\sqrt{N}} \sum_{k\in\bbZ} x_k \varphi\left( \frac{x}{N} -k
\right)\\
&=& \sum_{k\in\bbZ} \sum_{j\in\bbZ} x_k c_j \varphi (x-Nk -j) \\
&=& \sum_{l\in\bbZ} \left( \sum_{k\in\bbZ} x_k c_{l-Nk}\right)
\varphi(x-l) = \big( W_\varphi S f \big) (x), 
\end{eqnarray*}
which follows from the identity $(5)$. This is the desired identity
claimed in the proposition. Introducing the scaling operator of $(1)$, the
identity $(7)$ takes the equivalent form $UW_\varphi = W_\varphi S$, and
we say that $W_\varphi$ intertwines the two operators $U$ and $S$. It is
called the wave operator, and it plays a central role in computational
harmonic analysis. 
\end{proof}

\begin{rem}\label{egremark}
The conclusion of the proposition, and the usefulness of the
$S_i$-operators, are relevent even if the scaling identity $(2)$ does not
have a solution in $L^2(\bbR)$. This is demonstrated for example in
\cite{Jor99a} and \cite{JP98} for the case $N=4$ with the following
scaling identity 
\begin{equation}
\varphi(x) = 2 \varphi(4x) + 2 \varphi (4x-2). 
\end{equation}
Rewriting $(8)$ as 
\begin{equation}
\int h(x)d  \varphi(x) = \frac{1}{2} \left( \int
h\left(\frac{x}{4}\right)d \varphi(x) + \int h\left( \frac{x+2}{4}\right)
d\,\varphi(x) \right),  
\end{equation}
with $h$ continuous, we see that $(8)$ has a unique probability measure
$d \varphi$ as its solution. It has Hausdorff dimension $HD[d\,\varphi] =
\frac{1}{2}$, so it is a singular measure. 
Moreover, $(8)$ does not have a solution $\varphi$ in $L^2(\bbR)\setminus
\{0\}$. 
\end{rem}

As discussed in the previous section, 
the scaling function $\varphi\in L^{2}\left( \mathbb{R}\right) $ for an
orthogonal wavelet, if it
exists,
may be
obtained as a solution to the equation
\begin{equation}
\varphi\left( x\right) 
=\sqrt{N}\sum_{k\in\mathbb{Z}}c_{k}\varphi\left( Nx-k\right) ,
\end{equation}
but existence of $L^{2}\left( \mathbb{R}\right) $ solutions
requires special assumptions on $m$, or equivalently
$\left( c_{k}\right) $, from (\ref{FourierSeries}),
to which we now turn.
In particular, we now obtain the equivalences which lead to 
wavelet representations.

\begin{lem}
\label{LemS:wavelets.1}We have the following
three equivalences.
\begin{enumerate}
\item \label{LemS:wavelets.1(1)} $m\in L^{\infty}\left(\mathbb{T}\right)$
if and only if $S$ is bounded.
\item \label{LemS:wavelets.1(2)} $\displaystyle\frac{1}{N}\sum_{\substack{w\in\mathbb{T}\\w^{N}=z}}
\left|m\left(w\right)\right|^{2}=1\;\mathrm{a.a.}\;z\in\mathbb{T}$
if and only if $S$ is isometric.
\item \label{LemS:wavelets.1(3)} If $m_{1},m_{2}\in L^{\infty}$ are given,
then the corresponding operators $S_{1}$ and $S_{2}$
have orthogonal ranges if and only if for $\mathrm{a.a.}\;z\in\mathbb{T}$,
\[
\frac{1}{N}\sum_{w^{N}=z}
\overline{m_{1}\left(w\right)}\,m_{2}\left(w\right)=0.
\]
\end{enumerate}
\end{lem}

\begin{rem}
\label{RemS:wavelets.3}
$S^*S$ is a multiplication operator by
$\frac1N\sum_{\substack{w\in\mathbb{T}\\w^N=z}}\left|m\left(w\right)\right|^2$
but $SS^*$ is not; in fact $S$ is not normal, nor even hyponormal.
\end{rem}

\begin{proof}[Proof of Lemma \textup{\ref{LemS:wavelets.1}}]
The equivalence in (\ref{LemS:wavelets.1(1)}) is clear. For
(\ref{LemS:wavelets.1(2)}), observe that for
$f\in L^{2}\left(\mathbb{T}\right)$,
\begin{align*}
\left(S^{*}Sf\right)\left(z\right)&=\frac{1}{N}\sum_{w^{N}=z}
\overline{m\left(w\right)}\,\left(Sf\right)\left(w\right)\\
&=\frac1N\sum_{w^N=z}
\overline{m\left(w\right)}\,m\left(w\right)f\left(w^N\right)\\
&=\frac{1}{N}\sum_{w^{N}=z}
\overline{m\left(w\right)}\,m\left(w\right)f\left(z\right)\\
&=\frac1N\sum_{w^N=z}\left|m\left(w\right)\right|^2f\left(z\right).
\end{align*}
In particular, $S^{*}S=I$ precisely when the condition on
$m\left(z\right)$ in (\ref{LemS:wavelets.1(2)}) is satisfied.
Finally, if $S_{1}$ and $S_{2}$ are
given by $m_{1}$ and $m_{2}$, a similar calculation shows that
\[
\left(S_{1}^{*}S_{2}f\right)\left(z\right)=\frac{1}{N}\sum_{w^{N}=z}
\overline{m_{1}\left(w\right)}\,m_{2}\left(w\right)f\left(z\right),
\]
whence the equivalence (\ref{LemS:wavelets.1(3)}) becomes apparent.
\end{proof}

Before continuing, let us set aside the basic definitions from wavelet
theory which we need \cite{Dau92}. 

\begin{defn}
By a \emph{wavelet of scale $N$} we mean a finite set of functions
$\psi_i$, $i=1,\dots,N-1$, in $L^2\left(\mathbb{R}\right)$
such that the family
\[
\psi_{i,j,k}\left(x\right):=N^{\frac j2}\psi_{i}\left(N^{j}x-k\right),
\qquad j,k\in\mathbb{Z},
\]
satisfies
\[
\ip{f}{f} =
\int_{\mathbb{R}}
\left|f\left(x\right)\right|^2\,dx=\sum_{i,j,k}
\left|\ip{f}{\psi_{i,j,k}}_{L^2\left(\mathbb{R}\right)}\right|^2
\]
for all $f\in L^2\left(\mathbb{R}\right)$. 

It is an {\it orthogonal
wavelet} when the family $\psi_{i,j,k}$ forms an orthonormal basis for
$L^2\left(\mathbb{R}\right)$, equivalently, when every $||\psi_{i,j,k}||
=1$. For such wavelets 
there is a $1$--$1$ and explicit correspondence between
the family $\left( \psi_{i}\right) _{i=1}^{N-1}$ together with the
associated scaling function $\phi$, and
systems of so called {\it wavelet filter functions} $\left( m_{i}\right)
_{i=0}^{N-1}$ which are characterized by condition $(ii)$ of
Theorem~\ref{ProS:wavelets.2} (see \cite{BJendII,BJwave,Jorgen_min}). 

More generally, a {\it biorthogonal wavelet} consists of two families $\{
\psi_i \}$ and $\{\td{\psi}_i\}$ of $N-1$ functions in $L^2 \left(
\mathbb{R}\right)$ such that
\[
\ip{f}{g} = \sum_{i,j} \ol{\big<  \psi_{i,j,k}\,\big|\, f \big>} 
\big< \td{\psi}_{i,j,k}\,\big|\, g\big> 
\qfor f,g \in L^2\left(\mathbb{R}\right).
\]
These wavelets also have associated filter functions $\{ m_i \}$ and $\{
\td{m}_i\}$ which satisfy the condition specified in $(ii)$ of
Theorem~\ref{ProS:wavelets.3}.
\end{defn}

As a first consequence of the previous lemma we obtain the well-known
method (for instance see \cite{BJwave,Jorgen_min}) of generating
Cuntz-algebra representations from orthogonal wavelets. The {\it Cuntz
algebra} $\O_N$ is the universal $\ca$-algebra generated by the relations
in $(i)$ of Theorem~\ref{ProS:wavelets.2}. It has been studied extensively
by operator algebraists since the work \cite{Cun}. 

\begin{thm}
\label{ProS:wavelets.2}The following
three conditions are equivalent when the functions
$m_{0},\dots,m_{N-1}\in L^{\infty}\left(\mathbb{T}\right)$
are given
and operators $S_{0},\dots,S_{N-1}$ are defined by
$S_{i}f\left(z\right)=m_{i}\left(z\right)f\left(z^{N}\right)$.
\begin{enumerate}
\item \label{ProS:wavelets.2(1)}
$\left\{
\renewcommand{\arraystretch}{1.25}
\begin{array}{l}
\displaystyle S_{i}^{*}S_{j}=\delta_{i,j}I,\\
\displaystyle \sum_{i=0}^{N-1}S_{i}S_{i}^{*}=I.
\end{array}
\right.$
\item \label{ProS:wavelets.2(2)} 
The functions $m_{0},\dots,m_{N-1}$
on the torus $\mathbb{T}$
are the filter functions for an orthogonal wavelet of scale $N$.
In other words,
the matrix
\[
\frac{1}{\sqrt{N}}\left(m_{k}\left(e^{i\frac{2\pi l}{N}}z\right)\right)
_{k,l=0}^{N-1}
\]
is in
$\mathrm{U}_{N}\left(\mathbb{C}\right)
_{\mathstrut}
\;\mathrm{a.a.}\;z\in\mathbb{T}$.
\item \label{ProS:wavelets.2(3)}
$\displaystyle A_{k,l}\left(z\right)=
\frac{1}{N}\sum_{w^{N}=z}w^{-l}m_{k}\left(w\right)$
are the entries of a loop
$\mathbb{T}\rightarrow\mathrm{U}_{N}\left(\mathbb{C}\right)$,
i.e., a matrix function
\[
A\left(z\right)=\left(A_{k,l}\left(z\right)\right)\in
\mathrm{U}_{N}\left(\mathbb{C}\right)\qquad\mathrm{a.a.}\;z\in\mathbb{T}.
\]
\end{enumerate}
\end{thm}

\begin{proof}
The two Cuntz identities in (\ref{ProS:wavelets.2(1)}) correspond to the
orthonormality of the rows and columns in the
matrices of (\ref{ProS:wavelets.2(2)}). Indeed, the previous lemma shows
that the $S_{i}$ being isometries with pairwise
orthogonal ranges is the same as the rows being
orthonormal. On the other hand, for $f\in L^{2}\left(\mathbb{T}\right)$ we have
\begin{align*}
&\ip{\sum_{i=0}^{N-1}S_{i}S_{i}^{*}f}{f}
=\sum_{i=0}^{N-1}\ip{S_{i}^{*}f}{S_{i}^{*}f}\\
&\qquad=\sum_{i=0}^{N-1}\frac{1}{N^{2}}\int_{\mathbb{T}}
\sum_{w^{N}=z=w^{\prime\,N}}
\overline{m_{i}\left(w\right)}\,m_{i}\left(w'\right)
f\left(w\right)\,\overline{f\left(w'\right)}\,d\mu\left(z\right)\\
&\qquad=\frac{1}{N}\int_{\mathbb{T}}
\sum_{w^{N}=z=w^{\prime\,N}}
\underset{=\delta_{w,w'}}{\underbrace{
\left(\frac{1}{N}\sum_{i=0}^{N-1}
\overline{m_{i}\left(w\right)}\,m_{i}\left(w'\right)\right)
}}
f\left(w\right)\,\overline{f\left(w'\right)}\,d\mu\left(z\right).
\end{align*}
However, we can write $\left\|f\right\|_{2}^{2}=\frac{1}{N}\int_{\mathbb{T}}
\sum_{w^{N}=z}\left|f\left(w\right)\right|^{2}\,d\mu\left(z\right)$. Hence
the identity $\sum_{i=0}^{N-1}S_{i}S_{i}^{*}=I$
is equivalent to the condition
\[
\frac{1}{N}\sum_{i=0}^{N-1}
\overline{m_{i}\left(w\right)}\,m_{i}\left(w'\right)=\delta_{w,w'},
\]
for $N$th roots $w$ and $w'$ of $\mathrm{a.a.}\;z\in\mathbb{T}$.
This is equivalent
to column orthonormality.
For the equivalence of
(\ref{ProS:wavelets.2(2)}) and (\ref{ProS:wavelets.2(3)}),
consider the following
calculation:
\begin{align*}
\sum_{k=0}^{N-1}A_{i,k}\left(z\right)\,\overline{A_{j,k}\left(z\right)}
&=\frac{1}{N^{2}}\sum_{k}\sum_{w^{N}=z=w^{\prime\,N}}
w^{-k}m_{i}\left(w\right)w^{\prime\,k}\,\overline{m_{j}\left(w^{\prime}\right)}\\
&=\frac{1}{N}\sum_{w,w^{\prime}}\left(\frac{1}{N}\sum_{k=0}^{N-1}
\left(w^{-1}w^{\prime}\right)^{k}\right)
m_{i}\left(w\right)\,\overline{m_{j}\left(w^{\prime}\right)}\\
&=\frac{1}{N}\sum_{w^{N}=z}m_{i}\left(w\right)
\,\overline{m_{j}\left(w\right)}.
\end{align*}
Thus the matrix $A\left(z\right)$ is unitary precisely when
the matrix in (\ref{ProS:wavelets.2(2)}) is unitary.
\end{proof}

There are a number of advantages obtained by using the
matrix approach given by the $A\left(z\right)$,
including the fact that the filter functions $m_{i}$ can be recovered from $A$
(see Section \ref{S:co}).
For our purposes, this approach is helpful when
considering co-invariant subspaces, and it provides
motivation for our Fock-space construction. Now let us
turn to the new result here: namely, every
biorthogonal wavelet yields 
operators on Hilbert space satisfying simple identities
which contain the Cuntz relations in the special case
of orthogonal wavelets. There is an analogous matrix approach as well.

\begin{thm}
\label{ProS:wavelets.3}
The following conditions are
equivalent when the functions
$m_{0},\dots,m_{N-1}$,
$\tm _{0},\dots,\tm _{N-1}$,
and the corresponding operators
$S_{i}f\left(z\right):=m_{i}\left(z\right)f\left(z^{N}\right)$,
$\tS _{i}f\left(z\right):=\tm _{i}\left(z\right)f\left(z^{N}\right)$
are given.
\begin{enumerate}
\item \label{ProS:wavelets.3(1)}
$\left\{
\renewcommand{\arraystretch}{1.25}
\begin{array}{l}
\displaystyle S_{i}^{*}\tS _{j}=\delta_{i,j}I,\\
\displaystyle \sum_{i=0}^{N-1}S_{i}\tS _{i}^{*}=I.
\end{array}
\right.$
\item \label{ProS:wavelets.3(2)}
The functions 
$m_{0},\dots,m_{N-1}$,
$\tm _{0},\dots,\tm _{N-1}$
are the filter functions for a 
biorthogonal wavelet of scale $N$. In other words,
the matrices 
\[
\displaystyle 
\frac{1}{\sqrt{N}}\left(m_{k}\left(e^{i\frac{2\pi l}{N}}z\right)\right)
_{k,l=0}^{N-1}
\mbox{\quad
and
\quad}
\displaystyle 
\frac{1}{\sqrt{N}}\left(\tm _{k}\left(e^{i\frac{2\pi l}{N}}z\right)\right)
_{k,l=0}^{N-1}
\]
belong to $\mathrm{GL}_{N}\left(\mathbb{C}\right)$
and the adjoint of one is the inverse of the other for
$\mathrm{a.a.}\;z\in\mathbb{T}$.
\item \label{ProS:wavelets.3(3)}
The two matrix functions $A$ and $\tA $ with entries
\begin{align*}
A_{k,l}\left(z\right)
&=\frac{1}{N}\sum_{w^{N}=z}w^{-l}m_{k}\left(w\right),\\
\tA _{k,l}\left(z\right)
&=\frac{1}{N}\sum_{w^{N}=z}w^{-l}\tm _{k}\left(w\right)
\end{align*}
satisfy
\[
\sum_{k=0}^{N-1}
\overline{A_{k,i}\left(z\right)}\,
\tA _{k,j}\left(z\right)=\delta_{i,j},\qquad \mathrm{a.a.}\;z\in\mathbb{T},
\]
i.e.,
\[
A^{*}\tA =I\text{\qquad pointwise, }\mathrm{a.a.}\;z\in\mathbb{T},
\]
or
\[
\tA =A^{*\,-1},
\]
where
the function
$z\rightarrow A^{*}\left(z\right)$ denotes the adjoint matrix function
mapping
$\mathbb{T}\rightarrow\mathrm{GL}_{N}\left(\mathbb{C}\right)$.
\end{enumerate}
\end{thm}

\begin{proof}
With small adjustments we can follow the lines
of the previous proof. From the lemma, it
follows that the condition $S_{i}^{*}\tS _{j}=\delta_{i,j}I$ is
equivalent to the identity
\[
\frac{1}{N}\sum_{w^{N}=z}
\overline{m_{i}\left(w\right)}\,\tm _{j}\left(w\right)=\delta_{i,j},
\]
for $0\leq i,j\leq N-1$ and $\mathrm{a.a.}\;z\in\mathbb{T}$, while the identity
$\sum_{i=0}^{N-1}S_{i}\tS _{i}^{*}=I$ is a restatement of
\[
\frac{1}{N}\sum_{i=0}^{N-1}
\overline{m_{i}\left(w\right)}\,\tm _{i}\left(w'\right)=\delta_{w,w'},
\]
for $N$th roots $w$ and $w'$ of $\mathrm{a.a.}\;z\in\mathbb{T}$.
Thus the first two
conditions are equivalent. Finally, we can see
conditions
(\ref{ProS:wavelets.3(2)}) and (\ref{ProS:wavelets.3(3)})
are equivalent by following
the computation in the previous proof with $\tA _{i,k}\left(z\right)$
replacing $A_{i,k}\left(z\right)$.
\end{proof}

\begin{eg}
The matrix functions $A : \bbT \rightarrow \mathrm{GL}_N(\bbC)$ of
Theorems~\ref{ProS:wavelets.2} and \ref{ProS:wavelets.3} might be constant
even though the filter functions $\{ m_i\}_{i=0}^{N-1}$ are
non-constant. If $N=2$ and 
\[
\varphi(x) = \left\{ \begin{array}{cl}
1 & \mbox{$0 \leq x < 1$} \\
0 & \mbox{other $x\in\bbR$}
\end{array}\right.
\]
is the scaling function for the Haar wavelet, then $m_0 (z) =
\frac{1}{\sqrt{2}} \left( 1+z \right)$ and 
$A(z) = \frac{1}{\sqrt{2}}  \begin{sbmatrix} 1 & 1 \\ 1 & -1
\end{sbmatrix}
\in \mathrm{U}_2(\bbC)$. 

For the example in Remark~\ref{egremark}, $N=4$, $m_0 (z) = 1+z^2$ and 
\[
A(z) = \left[ \begin{matrix}
1&0&1&0 \\
1&0&-1&0 \\
0&1&0&1 \\
0&1&0&-1 
\end{matrix} \right]
\]
is an admissible non-unitary matrix function, here mapping $\bbT
\rightarrow \mathrm{GL}_4 (\bbC)$. This yields a biorthogonal system,
where the biorthogonality is encoded in the duality between measures and
continuous functions on $[0,1]\subset \bbR$. To get scaling functions
belonging to $L^2(\bbR)$, $m_0$ must satisfy $\frac{1}{N} \sum_{w^N = z} |
m_0 (w) |^2 \leq 1 $ \cite{BrJo01}; and one checks that $m_0 (z) = 1 +
z^2$ violates this. In fact, $\frac{1}{4} \sum_{w^4 =z} |m_0(w)|^2 = 2$. 

Our convention for the filter function $m_0 (z) = \sum_k c_k z^k$ is that
the coefficients $(c_k)$ are the masking numbers in $(10)$. Hence we
assume that $m_0(1) =\sqrt{N}$ where $N$ is the scaling number from
$(10)$. Introducing the $2\pi$-periodic variant of $m_0$, i.e.,
\[
m_0(t) := m_0\left(e^{-it}\right),
\]
we get the product formula for the Fourier transform 
\[
\hat{\varphi} (t) = \int_{\bbR} e^{-itx} d\varphi (x)
\]
from $(10)$ in the form 
\begin{equation}
\hat{\varphi}(t) = \prod_{j=1}^\infty \frac{m_0\left( t / N^j
\right)}{\sqrt{N}}, \qfor t\in\bbR.
\end{equation} 
This works even if $d\varphi$ is only a tempered distribution. If 
\[
\frac{1}{N} \sum_{k=0}^{N-1} \left| m_0\big( t+ \frac{k2\pi}{N}
\big)\right|^2 \leq 1,
\]
it follows that the infinite product in $(11)$ is in $L^2(\bbR)$. 
\end{eg}

\begin{defn}
We refer to the relations in $(i)$ of Theorem~\ref{ProS:wavelets.3} as the
{\it biorthogonal relations}. Further, any system $\{S_i\}$ satisfying the
relations in $(i)$ of Theorem~\ref{ProS:wavelets.2} clearly determines a
representation of $\O_N$, what we call an {\it orthogonal wavelet
representation}. Similarly, we will refer to a system $\{S_i, \td{S_j}\}$
satisfying the biorthogonal relations as a {\it biorthogonal wavelet
representation}.  
\end{defn}

\begin{rem}
\label{RemAfterProS:wavelets.3}
There are a number of properties which can
be derived from the biorthogonal relations. For instance, a simple
matrix argument shows there are no finite-dimensional  representations of
them. Further,
there are some nice reductions which can be made
on words in the generators
$\left\{S_{i},S_{i}^{*},\tS _{i},\tS _{i}^{*}\right\}$.
However,
overall there is little we can say about operators which satisfy these
relations in full generality. Fortunately, the biorthogonal 
wavelet representations of the biorthogonal relations
have some additional properties. As we shall see in the
next section, they have tractable finite-dimensional co-invariant cyclic
subspaces, and the generators satisfy other helpful
relations.
\end{rem}

The following simple lemma (see, e.g., \cite{BJendII}) will
help clarify the discussion below:

\begin{lem}
\label{LemS:wavelets.4}Let $R:=R_{N}$ be the average operator
\[
Rf\left(z\right):=\frac{1}{N}\sum_{\substack{w\in\mathbb{T}\\w^{N}=z}}
f\left(w\right).
\]
It is contractive in $L^{\infty}\left(\mathbb{T}\right)$,
and co-isometric in $L^{2}\left(\mathbb{T}\right)$.
Setting $S_{0}f\left(z\right)=f\left(z^{N}\right)$, we
have $R=S_{0}^{*}$, where the $*$
refers to the adjoint operation
relative to $L^{2}\left(\mathbb{T}\right)$, and
\[
\ker\left(R\right)^{\perp}=S_{0}L^{2}\left(\mathbb{T}\right).
\]
\end{lem}

\begin{proof}Let
$e_{k}\left(z\right):=z^{k}$, $z\in\mathbb{T}$, $k\in\mathbb{Z}$.
This is the standard Fourier
basis in $L^{2}\left(\mathbb{T}\right)$. Using duality
for the finite cyclic group
$\mathbb{Z}\diagup N\mathbb{Z}\cong\left\{0,1,\dots,N-1\right\}$, we then
get
\[
Re_{k}=\left\{
\renewcommand{\arraystretch}{1.25}
\begin{array}{ll}
e_{k/N}&\text{ if }k\equiv0\bmod{N},\\
0&\text{ if }k\not\equiv0\bmod{N}.
\end{array}
\right.
\]
The remaining details are left to the reader; or see \cite{BJendII}.
\end{proof}

Recall that if $L^2(\bbT)$ is a module over $\A$, then functions $\{ m_i
\} \subseteq L^2(\bbT)$ form a {\it module basis} for this module if given
$f\in L^2(\bbT)$ there is a unique expansion $f = \sum_i m_i a_i$ with
$a_i \in\A$. 

\begin{cor}
\label{CorS:wavelets.5}$\mathcal{A}_{1}
:=S_{0}\left(L^{\infty}\left(\mathbb{T}\right)\right)$ is a subalgebra of
$L^{\infty}\left(\mathbb{T}\right)$, and $L^{2}\left(\mathbb{T}\right)$ is a
module over $\mathcal{A}_{1}$ of module dimension $N$.
In fact, the functions
$m_{0},\dots,m_{N-1}$ form a module
basis for $L^{2}\left(\mathbb{T}\right)$ over $\mathcal{A}_{1}$
if and only if they satisfy condition \textup{(\ref{ProS:wavelets.2(2)})}
in
Theorem \textup{\ref{ProS:wavelets.2}}, i.e., if and only if
\[
R\left(\bar{m}_{i}m_{j}\right)=\delta_{i,j}1.
\]
\end{cor}

\begin{proof}
Since
\begin{align*}
R\left(\bar{m}_{i}m_{j}\right)\left(z\right)
&=\frac{1}{N}\sum_{w^{N}=z}\bar{m}_{i}\left(w\right)m_{j}\left(w\right)\\
&=\frac{1}{N}\sum_{k=0}^{N-1}
\bar{m}_{i}\left(e^{i\frac{2\pi k}{N}}z^{\frac{1}{N}}\right)
m_{j}\left(e^{i\frac{2\pi k}{N}}z^{\frac{1}{N}}\right)
\end{align*}
where $z^{\frac{1}{N}}$ is the principal branch
of the $N$'th root, the orthogonality
relations are clear. By
(\ref{ProS:wavelets.2(1)}) $\Leftrightarrow$ (\ref{ProS:wavelets.2(2)})
in
Theorem \ref{ProS:wavelets.2}, we have
\[
L^{2}\left(\mathbb{T}\right)\ni f\left(z\right)
=\sum_{i=0}^{N-1}S_{i}S_{i}^{*}f\left(z\right)
=\sum_{i=0}^{N-1}m_{i}\left(z\right)\left(S_{i}^{*}f\right)\left(z^{N}\right)
\in\sum_{i=0}^{N-1}m_{i}\mathcal{A}_{1}
\]
which shows that the orthogonality
relations make $m_{0},\dots,m_{N-1}$
a module basis.
\end{proof}

\begin{cor}
\label{CorS:wavelets.6}$\mathcal{A}_{k}
:=S_{0}^{k}\left(L^{\infty}\left(\mathbb{T}\right)\right)$ is a subalgebra of
$\mathcal{A}_{k-1}$ of module
dimension $N$, and the
module dimension of $L^{2}\left(\mathbb{T}\right)$ over
$\mathcal{A}_{k}$ is $N^{k}$. The corresponding module
basis is
\[
b_{i_{1},i_{2},\dots,i_{k}}
:=m_{i_{1}}\left(z\right)m_{i_{2}}\left(z^{N}\right)\cdots
m_{i_{k}}\left(z^{N^{k-1}}\right).
\]
\end{cor}

\begin{proof}
Every $f\in L^{\infty}\left(\mathbb{T}\right)$
satisfies
\[
f=\sum_{i_{1},\dots,i_{k}}S_{i_{1}}\cdots S_{i_{k}}
S_{i_{k}}^{*}\cdots S_{i_{1}}^{*}f.
\]
Setting $f_{i_{1},\dots,i_{k}}:=S_{i_{k}}^{*}\cdots S_{i_{1}}^{*}f$,
we get
\[
f\left(z\right)=\sum_{i_{1},\dots,i_{k}}b_{i_{1},\dots,i_{k}}\left(z\right)
f_{i_{1},\dots,i_{k}}\left(z^{N^{k}}\right)
\in\sum_{i_{1},\dots,i_{k}}b_{i_{1},\dots,i_{k}}\mathcal{A}_{k}.\settowidth
{\qedskip}{$\displaystyle f\left(z\right)
=\sum_{i_{1},\dots,i_{k}}b_{i_{1},\dots,i_{k}}\left(z\right)
f_{i_{1},\dots,i_{k}}\left(z^{N^{k}}\right)
\in\sum_{i_{1},\dots,i_{k}}b_{i_{1},\dots,i_{k}}\mathcal{A}_{k}.$}\addtolength
{\qedskip}{-\textwidth}\rlap{\hbox to-0.5\qedskip{\hfil\qed}}
\]
\renewcommand{\qed}{}
\end{proof}

\section{Co-invariant subspaces}\label{S:co}


{}From the theory of wavelets \cite{Dau92},
a compactly supported biorthogonal wavelet of scale $N$ is determined by
scaling functions $\varphi$ and $\tilde{\varphi}$ which generate the associated
$2N-2$ wavelet functions. Further, the functions $\varphi,\tilde{\varphi}$ are
supported on the interval
$\left[-Ng+1,Ng-1\right]$, where $g$ is the \emph{genus}
of the wavelet.
(We point out that the recent paper \cite{BrJo01} examines
spaces of all such scaling functions.)
In this case,
the corresponding filter functions $m_i$ and $\tm _i$ are
Fourier
polynomials of degree $Ng-1$,
i.e., of the form $\sum_{k=-Ng+1}^{Ng-1}a_{k}z^{k}$.
The numbers $a_{k}$ are the wavelet
masking coefficients, i.e.,
$\varphi\left( x\right) =\sum_{k}a_{k}\varphi\left( Nx-k\right) $.
If $a_{k}=0$ unless $0\leq k\leq Ng-1$,
then $\varphi$ is supported in $\left\lbrack 0,Ng-1\right\rbrack $.
One of the advantages of using the matrix perspective
given by the invertible loops
is that this degree is considerably reduced for the functions
$A_{i,j}$, $\tA _{i,j}$. We begin by observing this fact, together with
the method of recovering the filter functions (hence the
$S$, $\tS $ system) from the invertible loops.

\begin{lem}
\label{LemS:co.1}
For $0\leq i\leq N-1$, the filters $m_i\left(z\right)$,
$\tm _i\left(z\right)$ are
obtained from $A$, $\tA $ by
\[
m_i\left(z\right)=\sum_{j=0}^{N-1}A_{i,j}\left(z^N\right)z^j
\text{\quad and\quad}
\tm _i\left(z\right)=\sum_{j=0}^{N-1}\tA _{i,j}\left(z^N\right)z^j.
\]
Further, if each $m_i$, $\tm _i$ is a polynomial of degree
at most $Ng-1$, then the $A_{i,j}$, $\tA _{i,j}$ have degree at
most $g-1$.
\end{lem}

\begin{proof}
The first claim simply follows from the computation
\begin{align*}
\sum_jA_{i,j}\left(z^N\right)z^j
&=\sum_j\left(\frac1N\sum_{w^N=z^N}m_i\left(w\right)w^{-j}\right)z^j\\
&=\sum_{w^N=z^N}m_i\left(w\right)\delta_{w,z}=m_i\left(z\right).
\end{align*}
The same computation works for $\tm _i$ and $\tA _{i,j}$. To verify
the second claim, suppose
$m_i\left(z\right)=\sum_{k=0}^{Ng-1}a_{k}^{\left( i\right) }z^k$. Then we have
\begin{align*}
A_{i,j}\left(z\right)
&=\frac1N\sum_{w^N=z}m_i\left(w\right)w^{-j}\\
&=\sum_{\left| k\right| \leq Ng-1}
a_{k}^{\left( i\right) }\left(\frac1N\sum_{w^N=z}w^{k-j}\right)\\
&=\sum_{\left| k\right| \leq Ng-1}
a_{k}^{\left( i\right) }z^{\frac{k-j}{N}}
\delta_{k-j\left(\operatorname*{mod}N\right),0}\,.
\end{align*}
However, the quantity $\frac{k-j}{N}$ is bounded above by $g-1$,
as required.
\end{proof}

The biorthogonal wavelet representations are rather specialized in that
they have tractable finite-dimensional, co-invariant subspaces
which are also \emph{doubly-cyclic}. We shall discuss the significance
of this fact in Remark 3.4. This comes as a direct consequence
of the following result. The key technical device is that, as for the
orthogonal wavelet representations \cite{BJwave}, the actions of the
adjoint operators on Fourier basis vectors can be
computed directly.

\begin{lem}
\label{LemS:co.2}
Let $S=\left(S_0,\dots,S_{N-1}\right)$,
$\tS =\left(\tS _0,\dots,\tS _{N-1}\right)$ be a compactly
supported biorthogonal wavelet representation of genus $g$.
Let $\left\{e_n:n\in\mathbb{Z}\right\}$ be the
basis for $L^2\left(\mathbb{T}\right)$ given by
$e_n\left(z\right)=z^n$ and let
\[
\mathcal{K}=\operatorname*{span}\left\{e_0,e_{-1},\dots,e_{-Ng+1}\right\}
\bigwedge\operatorname*{span}_{i,j}\left\{\overline{A_{i,j}\left( z\right)
}
\,z^{r}\mathrel{;}r\leq 0\right\}.
\]
Then we have
\begin{equation}
S_i^*\mathcal{K}\subseteq\mathcal{K}
\text{\quad and\quad}
\tS _i^*\mathcal{K}\subseteq\mathcal{K}
\text{\quad for }
0\leq i\leq N-1.
\label{CoInvariance}
\end{equation}
\textup{(}Properties \textup{(\ref{CoInvariance})}
are called \emph{co-invariance}.\textup{)}
Further, for all $n\in\mathbb{Z}$ there is a $K\geq1$ such that
\[
S_{i_1}^*\cdots S_{i_k}^*e_n\in\mathcal{K}
\text{\quad and\quad}
\tS _{i_1}^*\cdots \tS _{i_k}^*e_n\in\mathcal{K}
\]
for all $k\geq K$ and all indices $0\leq i_1,\dots,i_k\leq N-1$.
\end{lem}

\begin{proof}
We first compute the action of $S_i^*$ on a typical basis
vector $e_n$. From the previous lemma we have
\begin{align*}
S_i^*e_n\left(z\right)
&=\frac1N\sum_{w^N=z}\overline{m_i\left(w\right)}\,w^n\\
&=\sum_{j=0}^{N-1}\overline{A_{i,j}\left(z\right)}\,
\left(\frac1N\sum_{w^N=z}w^{n-j}\right)\\
&=\sum_{j=0}^{N-1}\overline{A_{i,j}\left(z\right)}\,
z^{\frac{n-j}{N}}\delta_{n-j\left(\operatorname*{mod}N\right),0}\,.
\end{align*}
Of course, only one term in this sum is nonzero.
Recall that the $A_{i,j}\left(z\right)$ are of degree at most $g-1$. Thus,
if $0\leq n\leq Ng-1$, it follows that $S_i^*e_{-n}$ is the complex
conjugate of a polynomial of degree at most $Ng-1$.
This says precisely that $S_i^*e_{-n}\in\mathcal{K}$.

It remains to check that the adjoints `pull back' basis
vectors into $\mathcal{K}$. The pattern becomes clear after two
applications of adjoints. Let
$A_{p,q}\left(z\right)=\sum_{l=0}^{g-1}A_{l}^{\left(p,q\right)}z^l$
for $0\leq p,q\leq N-1$.
Let $n\in\mathbb{Z}$ and set $j_0\equiv n\left(\operatorname*{mod}N\right)$
with $0\leq j_0\leq N-1$.
Then we have
\begin{align*}
\overline{S_p^*S_i^*e_n}\,\left(z\right)
&=\frac1N\sum_{w^N=z}m_p\left(w\right)A_{i,j_0}\left(w\right)
w^{-\frac{n-j_0}{N}}\\
&=\frac1N\sum_{w}
\left(\sum_qA_{p,q}\left(z\right)w^q\right)
\left(\sum_lA_{l}^{\left(i,j_0\right)}w^l\right)
w^{-\frac{n-j_0}{N}}\\
&=\sum_{q,l}A_{p,q}\left(z\right)A_{l}^{\left(i,j_0\right)}
\left(\frac1N\sum_{w^N=z}w^{q+l-\frac{n-j_0}{N}}\right)\\
&=\sum_{q,l}A_{l}^{\left(i,j_0\right)}
\delta_{q+l-\frac{n-j_0}{N}\left(\operatorname*{mod}N\right),0}
A_{p,q}\left(z\right)z^{\frac{q+l-\frac{n-j_0}{N}}{N}}.
\end{align*}
{}From this computation we see that vectors $S_{i_1}^*\cdots S_{i_k}^*e_n$
belong to the span of vectors of the form
$\overline{A_{p,q}\left(z\right)}\,z^r$,
where the absolute value of the powers $r$
decreases steadily as $k$ increases. It follows that
$S_{i_1}^*\cdots S_{i_k}^*e_n$ will belong to $\mathcal{K}$ when $k$ is large
enough, and co-invariance means it will stay there.

We have carried out the analysis on the $S_i^*$, but
the same proof works for the $\tS _i^*$. From Lemma \ref{LemActionAdjoint}
they have analogous formulae,
and from the discussion at the start of this section
the same compact support yields the same
summation limits throughout.
\end{proof}

\begin{thm}
\label{ThmS:co.3}
Let $S=\left(S_0,\dots,S_{N-1}\right)$,
$\tS =\left(\tS _0,\dots,\tS _{N-1}\right)$ be
a compactly supported biorthogonal wavelet representation
on $\H=L^2\left(\mathbb{T}\right)$. Then there is a finite-dimensional subspace
$\mathcal{K}$ which is co-invariant and doubly-cyclic for the representation.
In other words,
\begin{enumerate}
\item \label{ThmS:co.3(1)}
$S_i^*\mathcal{K}\subseteq\mathcal{K}$ and
$\vphantom{\tilde{S}_i^*}
\tS _i^*\mathcal{K}\subseteq\mathcal{K}$ for $0\leq i\leq N-1$,
\item \label{ThmS:co.3(2)}
$\displaystyle
\bigvee_{\substack{i_1,\dots,i_k\\k\geq1}}^{{}}
S_{i_1}\cdots S_{i_k}\mathcal{K}
=\H
=\bigvee_{\substack{i_1,\dots,i_k\\k\geq1}}^{{}}
\tS _{i_1}\cdots \tS _{i_k}\mathcal{K}$.
\end{enumerate}
\end{thm}

\begin{proof}
The subspace $\mathcal{K}$
from Lemma \ref{LemS:co.2}
provides the
candidate.
The only thing left to show is cyclicity. But
from the biorthogonal relations, for $k\geq1$
we have
\begin{equation}
\sum_{0\leq i_1,\dots,i_k\leq N-1}S_{i_1}\cdots S_{i_k}
\tS _{i_k}^*\cdots\tS _{i_1}^*=I.
\label{eqSumRule}
\end{equation}
In particular,
when this identity is applied to Fourier basis vectors,
the previous lemma yields (\ref{ThmS:co.3(2)}) for $S$.
Specifically, let $n\in\mathbb{Z}$. Using
then Lemma \ref{LemS:co.2}, we may pick
some $k\in\mathbb{N}$, depending on $n$,
such that
\[
\tS _{i_k}^{\ast}\cdots \tS _{i_1}^{\ast}e_{n}\in\mathcal{K}
\text{\qquad for all }i_1,\dots,i_k.
\]
An application of (\ref{eqSumRule}) then
yields $e_{n}\in\bigvee_{i_1,\dots,i_k}
S_{i_1}\cdots S_{i_k}\mathcal{K}$.
The result follows since the closed
span of the monomials $e_{n}$, $n\in\mathbb{Z}$, is
$L^2\left(\mathbb{T}\right)$.
Finally, taking the adjoint of this identity and applying
the lemma again for $\tS $ completes the proof.
\end{proof}

\begin{rem}
\label{RemS:co.4}
There is an entire structure theory for representations of $\O_N$
which are found to have finite-dimensional, co-invariant
cyclic subspaces. Indeed, the recent paper \cite{DKS}
sets out an entire theory for decomposing such
representations into tractable classes of irreducible
subrepresentations. This paper was presented in the context
of dilation theory, but it was observed
in \cite{Jorgen_min} and \cite{Kri_cp} that the orthogonal wavelet
representations of
$\mathcal{O}_N$
form a subclass of these representations.

For such a representation, let $A_i$ be the compressions
of the isometries $S_i$ to a given finite-dimensional, co-invariant
cyclic subspace.
The crucial point in the analysis is that the finite-dimensional
\emph{minimal}
$A_i^*$-invariant subspaces
generate
the irreducible subspaces
for
the representation.
This came from an investigation into the completely positive
map $\Phi\left(X\right)=\sum_{i=0}^{N-1}A_iXA_i^*$ determined by the $A_i$.
Thus
decomposing these representations, acting on infinite-dimensional
space, amounts to computing for these finite-dimensional
`anchor' subspaces.
In fact, the paper \cite{Kri_cp} shows  these subspaces can
be obtained through a relatively simple analysis of the
map $\Phi$, without \emph{any} explicit reference to
the compressions $A_i$.

There are obvious analogues of this theory for the
biorthogonal setting, but one immediately confronts
serious issues.
The analogue of $\Phi$ would be a \emph{completely bounded}
unital map $\Phi\left(X\right)=\sum_{i=0}^{N-1}A_iX\tA _i^*$.
But in the orthogonal-$\O_N$-completely positive setting,
the key technical device is the unique dilation theory
which abounds: namely, what's known as the Frahzo--Bunce--Popescu
unique minimal isometric dilation of a row
contraction \cite{Bun,Fra1,Pop_diln},
which is really a special case of
Stinespring's unique
dilation of a completely positive map to a $C^*$-homomorphism
\textup{(}see \cite{Pa_cb_book}\textup{)}.
This allows us to go back and forth interchangeably between the $A_i$ and $S_i$,
as well as, respectively, the completely positive map and the
endomorphism determined by these operators.
In our more general setting dilations typically exist,
but they are \emph{not unique}. For instance, recall
from Paulsen's `off-diagonal technique' \cite{Pa_cb_book}
how completely bounded maps are dilated:
Every completely bounded map can be regarded as the off-diagonal
corner of a completely positive map. Stinespring's dilation
theorem gives a unique dilation of this map, which in turn
yields a completely bounded homomorphism which dilates
the completely bounded map.
The problem is that the way in which the map is regarded as an
off-diagonal corner is not unique \textup{(}in fact an application of
Arveson's matricial
Hahn-Banach Theorem is involved \cite{ArvI,ArvII}\textup{)}.
There is also the issue of irreducibility here. It is
not even clear what it should mean for an $S,\tS $ system
to be irreducible. 

Nonetheless, for the wavelet representations
of the biorthogonal relations
at least, Theorem \textup{\ref{ThmS:co.3}} shows that a weaker
spatial version of these dilation results is valid here.
In particular, the representations can be recovered spatially from the
compressions to particular finite-dimensional anchor
subspaces. The reader may notice that the computations
above can be strengthened to reduce the size of $\mathcal{K}$.
In fact, it appears that the analogue here of Section 8
from \cite{Jorgen_min} could be developed to obtain `minimal' subspaces
$\mathcal{L}$ and $\widetilde{\mathcal{L}}$ of $\mathcal{K}$
defined respectively by $A$ and $\tA $, 
which are 
co-invariant for $S$ \textup{(}respectively $\tS $\/\textup{)}
and cyclic for $\tS $ \textup{(}respectively $S$\/\textup{)}.
{}From \cite{DKS,Jorgen_min}, an orthogonal wavelet representation is
irreducible exactly when there is a unique such $\mathcal{L}$.
It would be interesting to know if there is an analogue
of this fact for the $\mathcal{L}$ and $\widetilde{\mathcal{L}}$ here.
\end{rem}

We finish this section by discovering a striking relationship between the
operators $\{S_i, \td{S}_j\}$ on the one hand, and the matrices
$A$ and $\tA $ on the other.

\begin{lem}
\label{LemS:wavelets.7}
Let $S=\left(S_0,\dots,S_{N-1}\right)$,
$\tS =\left(\tS _0,\dots,\tS _{N-1}\right)$ be
a biorthogonal wavelet representation
with invertible-loop matrix functions $A$, $\tA$.
Then for $f\in L^2\left(\mathbb{T}\right)$
and $0\leq i,j\leq N-1$
we have
\begin{enumerate}
\item \label{LemS:wavelets.7(1)}
$\displaystyle S_{i}^{*}S_jf\left(z\right)
=\left(AA^{*}\right)_{j,i}\left(z\right)f\left(z\right)$
and
\item \label{LemS:wavelets.7(2)}
$\displaystyle \vphantom{\tilde{S}_{i,j}^{*}}\tS_{i}^{*}\tS_jf\left(z\right)
=\left(AA^{*}\right)_{j,i}^{-1}\left(z\right)f\left(z\right)$.
\end{enumerate}
\end{lem}

\begin{proof}
For the $S_{i}$ we have the following computation:
\begin{align*}
S_{i}^{*}S_{j}f\left(z\right)
&=\frac{1}{N}\sum_{w^{N}=z}
\bar{m}_{i}\left(w\right)m_{j}\left(w\right)f\left(z\right)\\
&=\frac{1}{N}\sum_{w^{N}=z}\sum_{k,l}\bar{A}_{i,k}\left(z\right)
\bar{w}^{k}A_{j,l}\left(z\right)w^{l}f\left(z\right)\\
&=\sum_{k,l}\delta_{k,l}\bar{A}_{i,k}\left(z\right)
A_{j,l}\left(z\right)f\left(z\right)\\
&=\sum_{k}\bar{A}_{i,k}\left(z\right)
A_{j,k}\left(z\right)f\left(z\right)\\
&=\left(AA^{*}\right)_{j,i}\left(z\right)f\left(z\right).
\end{align*}
A similar calculation shows that
\[
\tS_{i}^{*}\tS_jf\left(z\right)
=\left(\tA \tA^{*}\right)_{j,i}\left(z\right)f\left(z\right)
=\left(AA^{*}\right)_{j,i}^{-1}\left(z\right)f\left(z\right).
\]
This completes the proof.
\end{proof}

\begin{rem}
\label{RemAfterLemS:wavelets.7}
This relationship provides us with the impetus for
our general Fock-space Hilbert space construction.
In particular,
the
$2N\times2N$ positive matrix $\mathcal{S}^{*}\mathcal{S}$, where
$\mathcal{S}=\left[S,\tS \right]$ and
$S=\left(S_{0},\dots,S_{N-1}\right)$,
$\tS =\left(\tS _{0},\dots,\tS _{N-1}\right)$ form a
biorthogonal wavelet representation, has the tractable form
\[
\mathcal{S}^{*}\mathcal{S}=
\begin{bmatrix}
AA^{*} & I_{N} \\
I_{N} & \left(AA^{*}\right)^{-1}
\end{bmatrix}
..
\]
Further,
this positive matrix
has
\textit{commuting}
entries since
the operators in the lemma are
\emph{multiplication} operators.

\end{rem}

\section{Fock space on positive matrices}\label{S:Fock}


There are now several Fock space constructions which appear in the
literature.
Typically, they allow certain identities to be
represented
by operators
on Hilbert space by way of  natural
left creation operators associated with the
underlying Fock space. See 
\cite{BS,JPS01,JSW95,JoWe94,KP,MN}
for some different perspectives. 
The purpose of this section is to introduce a new Fock space
construction which, we believe,
may provide the appropriate framework for
studying the biorthogonal wavelet representations discussed
above. In any event, we find this construction
to be interesting in its own right. To establish the nomenclature we
use for the next two sections, we begin by
reviewing the formulation of unrestricted
Fock space. In Example 4.12 below we point out how this motivating special
case fits into our construction.  
  
\begin{eg}
\label{Suggestive}
The \emph{full} (\emph{unrestricted}\/) \emph{Fock space}
over $\bbC^N$, where $N$ is a fixed positive integer
with $N\geq2$, is
the orthogonal direct sum of Hilbert spaces given by:
\[
\K = 
\left( \mkern6mu
\sideset{}{^{\smash{\oplus}}}{\sum}\limits_{k=-\infty}^{-1}
(\bbC^N)^{\otimes -k}  \right) \oplus \bbC
= \ldots \oplus (\bbC^N \otimes \bbC^N) \oplus (\bbC^N) \oplus
\bbC.
\]
The number {\bf 1} in the summand on the right
(giving the copy of $\mathbb{C}$)
is called the {\it
vacuum
vector} and is denoted by $\Omega$.
Let $\{ \xi_1, \ldots ,\xi_N\}$ be a fixed orthonormal basis for $\bbC^N$.
Then $\K$ is an infinite-dimensional Hilbert space with orthonormal basis
given by 
\[
\big\{ \xi_{i_1}\otimes \ldots \otimes \xi_{i_k} \,\,\big|\,\, 1 \leq
i_1,\ldots,i_k\leq
N,\, k\geq 1\big\} \cup \{\Omega\}. 
\]
We wish to think of the infinite direct
sum as extending from right to left.
This is non-standard, but we believe it is helpful in
understanding the action of the creation operators (see below),
and it allows us to introduce notation which is
less cumbersome. 

The {\it left creation operator} $L_i $ determined by $\xi_i$ on
$\K$ is defined by the actions: 
\[
\left\{
\renewcommand{\arraystretch}{1.25}
\begin{array}{l}
L_i ( \Omega)  = \xi_i \\
L_i (\eta_k \otimes \dots \otimes \eta_1 )
=  \xi_i \otimes \eta_k \otimes \dots
\otimes  \eta_1 ,
\end{array}
\right.
\]
for all $k \geq 1$ and $\eta_1, \dots, \eta_k \in \bbC^N$.
The adjoint of $L_i$ is the {\it annihilation operator}
determined by $\xi_i$, and it acts by:
\[
\left\{
\renewcommand{\arraystretch}{1.25}
\begin{array}{l}
L_i^* ( \Omega ) =  0 \\
L_i^* (\eta_1 ) = \ip{\eta_1}{\xi_i} \Omega \\
L_i^* (\eta_k \otimes \dots \otimes \eta_1 )  =  
\ip{\eta_k}{\xi_i} ( \eta_{k-1}
\otimes \dots \otimes \eta_1 ),
\end{array}
\right.
\]
for all $k \geq 2$ and $\eta_1, \dots, \eta_k \in \bbC^N$.
This terminology comes from theoretical physics where `creation' signifies
the creation of a new particle. 

There is another formulation of Fock space $\K$ which leads to a
notational simplification for us. 
Let $\bbF_N^+$ be the
unital free semigroup
on $N$ \emph{non-commuting} letters $\{ 1,2, \dots, N\}$ with
unit $e$.
Given $w$ in $\bbF_N^+$,
the positive integer
$\left|w\right|$ is the length of
the word $w$. The unit
$e$ corresponds to
the word of length
zero, or the empty
word.
Then one can also think of the Fock space $\K$ as
$\ell^2 (\bbF_N^+)$, where an orthonormal basis is given by
the vectors $\{ \xi_w : w\in \bbF_N^+\}$ corresponding to words.
Thus the vectors $\xi_{i_1} \otimes \ldots \otimes \xi_{i_k}$ are
identified with $\xi_w$ where
the product $w= i_1 \cdots i_k$ is in the free semigroup $\bbF_N^+$. Also, 
the vacuum vector is identified with $\xi_e$. We shall further simplify
notation by referring to the vector $\xi_w$ just by the word $w$. 
Hence the action of the creation operators is 
encapsulated in
the short statement
\[
L_i (w) = iw \qfor w \in \bbF_N^+, 
\]
where again we emphasize that the product $iw$ is in the free semigroup
$\bbF_N^+$. 
The actions of the annihilation operators are also easily described by   
$L_i^* (e ) = 0$, and $L_i^* (jw)
= w $ when
$i=j$ and 0 otherwise. 
These operators can, in fact,
be defined independent of basis (for $\xi\in\bbC^N$, an
operator  $L_\xi$ can be analogously defined).
We shall see this is also the case in our setting,
but it is convenient to fix a basis
for
the analysis.

It is not hard to see that $L= (L_1, \dots, L_N)$ forms an
$N$-tuple of isometries with pairwise
orthogonal
ranges, for which the closed span of the ranges of the isometries span the
orthogonal complement of the span of the vacuum vector.
Equivalently, since the $L_i L_i^*$
are the range projections, this says
\[
L_i^* L_j = \delta_{i,j} I \qfor 1\leq i,j \leq N
\qand \sumiN L_i L_i^* = I - P_{\Omega}. 
\]
These are the so called \emph{Cuntz--Toeplitz}
isometries, and the $\ca$-algebra they generate is
denoted by $\E_N$.
The ideal generated by the rank one projection $P_{\Omega}$ in $\E_N$
determines a copy of the compact operators,  and when it is factored out
the Cuntz algebra $\O_N$ is obtained.
Thus there is a tight relationship between $\O_N$
and the operators $L = (L_1, \ldots, L_N)$.  
\end{eg}

\begin{note}
\label{Confusion}
The reader will notice that in the previous example, and
for the next two sections, we have changed our notation
with $N$-tuples of operators from $\left\{0,1,\dots,N-1\right\}$ to
$\left\{1,2,\dots,N\right\}$. Unfortunately, this is the price to pay for
combining the two different perspectives. In wavelet analysis the
standard notation for multiresolution wavelet functions is
the former \textup{(}$0$ is for `low frequency'\/\textup{)},
while in the realm of theoretical physics and creation operators the
latter is necessary to portray the `creation' of new particles.
In any event, we hope this note will preempt any confusion.
\end{note}

The starting point for our general construction is
an extension result for  completely
positive maps. First,
let us recall the dichotomy between completely positive maps and positive
matrices given by Choi's Lemma \cite{Choi}. Let
$\{ e_{i,j} \}_{1\leq i,j \leq N}$ be
matrix units for the set of $N \times N$ complex matrices $\M_N$
corresponding to
a fixed orthonormal
basis $\left\{\xi_{1},\dots,\xi_{N}\right\}$ for
$\mathbb{C}^{N}$. Our construction
is independent of
basis (see Remark
\ref{RemBefore4.6}), but
for the sake of
brevity we shall
work with a fixed
basis.
The completely positive maps $\Phi \colon \M_N \rightarrow \bofh$
can be identified with the positive matrices
$P = [p_{i,j}] \in \M_N (\bofh)$, where the
correspondence is given by 
\[
P = \Phi^{(N)} \big( [e_{i,j}] \big) = \big[ \Phi (e_{i,j}) \big].
\]
We shall call $P= \big[ \Phi (e_{i,j}) \big]$ the
\emph{Choi matrix} associated with $\Phi$. Every
such completely positive map can be extended
in a natural way to the matrix algebras
$\M_{N^k}$.

\begin{lem}
\label{LemTildePhi}
Let $\Phi \colon \M_N \rightarrow \bofh$ be a completely positive map.
Then there is a unique
map $\td{\Phi} \colon \bigcup_{k\geq 1} \M_{N^k} \rightarrow \bofh$ such that
\[
\td{\Phi} ( a\otimes b) = \Phi(a) \td{\Phi} (b) 
\]
whenever $a\in\M_N$ and
$b\in \bigcup_{k\geq 1} \M_{N^k}$. In particular, for $a_1, \dots, a_k
\in \M_N$ we have $\td{\Phi}(a_1 \otimes \dots \otimes a_k)
= \Phi(a_1) \cdots \Phi(a_k)$. The
natural extension of $\Phi$ to $\uhfNI$ is not necessarily bounded.
\end{lem}

\begin{proof}
The definition of $\td{\Phi}$ is forced upon us by the conclusion.
In fact, $\td{\Phi}$ is
determined by a sequence of completely positive maps on the algebras
$\M_{N^k}$. The
matrices
$e_{i_1 i_1'} \otimes\dots\otimes e_{i_k i_k'}$, for $1 \leq i_j, i_j' \leq N$
and $1\leq j \leq k$,
form a set of matrix units for the $k$-fold tensor algebra
$\M_{N^k} \cong \M_N^{\otimes k}$.
(We use the standard
identification of
matrices
in $\M_{N^k}$ with tensors
found in such texts
as \cite{Pa_cb_book}.)
As above, let $P = [p_{i,j}] = \Phi^{(N)} \big( [e_{i,j}]
\big)$.
For $k\geq 1$, define maps $\Phi_k \colon \M_{N^k} \rightarrow \bofh$ by
\[
\Phi_k ( e_{i_1 i_1'} \otimes\dots\otimes e_{i_k i_k'})
= \Phi(e_{i_1 i_1'}) \cdots \Phi(e_{i_k
i_k'}) = p_{i_1 i_1'} \cdots p_{i_k i_k'}.
\]
Each of these maps is completely positive by Choi's Lemma since
\[
\Phi_k^{(N^k)} \big( [e_{i_1 i_1'} \otimes\dots\otimes e_{i_k i_k'}] \big)
=  [ p_{i_1 i_1'}
\cdots p_{i_k i_k'} ] 
\cong
P^{\otimes k} \geq 0,
\]
where the indices in the first two matrices satisfy $1\leq i_j,i_j' \leq N$
and $1\leq j \leq k$.
Thus for $a\in\M_{N^k}$,
define $\td{\Phi}(a) = \Phi_k (a)$. Then $\td{\Phi}\colon \bigcup_{k\geq 1}
\M_{N^k} \rightarrow \bofh$ is a  map which has the desired properties.
Uniqueness clearly follows from these properties. 

When $\bigcup_{k\geq 1} \M_{N^k}$
is regarded as an increasing union (given by unital
embeddings) which generates $\uhfNI$,
the natural extension of $\Phi$ to $\uhfNI$ will be
unbounded in general.
Indeed,
the identity in this algebra is obtained as a limit $I =
\lim_{k\rightarrow\infty} I_{N^k}$, and for $k \geq 1$,
$I_{N^k}\cong I_{N}^{\otimes k}$.
Hence we would have 
\[
\td{\Phi} (I) = \lim_{k\rightarrow\infty}
\td{\Phi} \left(I_{N^k}\right)
=
\lim_{k\rightarrow\infty} \Phi(I_N)^k,
\]
which may be unbounded if $\Phi$ is not completely contractive.
\end{proof}

We
are not concerned with the viability of an extension to
$\uhfNI$ since it is not necessary for the
Fock-space construction. The crucial point for us is that
completely positive maps
on $\M_{N}$
can be
extended to the
`pre-$\uhfNI$' algebras $\M_{N^k}$. We will let $\Phi$ denote the map
{\it and} its extension when there is no confusion. 

\begin{const}
\label{ConS:Fock.3}Let $\H$ be a Hilbert space.
Heuristically, our construction can be thought of as formally taking
the tensor product of unrestricted Fock space with $\H$,
then defining a `twisted' inner product
on the result by using a completely positive map from the
complex
matrices into $\bofh$ (or, if you like, a
positive matrix with entries in $\bofh$).

We define the
$N$\emph{-variable pre-Fock space over} $\H$
to be the vector space of finite sums
\[
\T_N (\H) = \left\{ \sum_{|w|\leq k} w\otimes h_w \Bigm| w\in\bbF_N^+,
 \; k\geq 1, \; h_w \in\H
\right\},
\]
where philosophically a vector $(i_1 \cdots i_k)\otimes h$, with $i_1,
\ldots, i_k \in \bbF_N^+$, corresponds to the vector 
$\xi_{i_1}\otimes\dots
\otimes \xi_{i_k} \otimes h$ in $\big( \bbC^N \big)^{\otimes k} \otimes
\H$. 
Let $\Phi \colon \M_N \rightarrow \bofh$ be a completely positive map
\textup{(}At this point we make no requirement
that $\Phi$ be completely bounded.\/\textup{)}.
Define a form
$\ip{\,\cdot\,}{\,\cdot\,}_\Phi \colon
\T_N (\H) \times \T_N (\H) \rightarrow \bbC$ in the following manner: For
$w$, $w'$ in $\bbF_N^+$ and $h$, $h'$ in $\H$, 
\begin{enumerate}
\item \label{ConS:Fock.3{1}}$\displaystyle\ip{e\otimes h}{e\otimes h'}_{\Phi}
=\ip{h}{h'}$;
\item \label{ConS:Fock.3(2)}if $|w| \neq |w'|$,
then $\displaystyle\ip{w\otimes h}{w' \otimes h'}_\Phi =  0$;
\item \label{ConS:Fock.3(3)}if $w=i_1\cdots i_k$ 
and $w' = i_1' \cdots i_k'$, then 
\[
\ip{w\otimes h}{w' \otimes h'}_\Phi = 
\ip{h}{\Phi( e_{i_1 i_1'} \otimes\dots\otimes e_{i_k i_k'}) h'}.
\]
\end{enumerate}
Then extend $\ip{\,\cdot\,}{\,\cdot\,}_\Phi$ to
$\T_N (\H) \times \T_N (\H) $ as linear in the first variable
and conjugate linear in the second.
\end{const}

\begin{thm}\label{possemidef}
The form $\ip{\,\cdot\,}{\,\cdot\,}_\Phi$ is
positive semi-definite on $\T_N (\H)$. 
\end{thm}

\begin{proof}
Let $x\in\T_N(\H)$ be a {\it finite} sum of the form
\[
x = \sum_{k \geq 0} \sum_{|w| = k} w\otimes h_w.
\]
As above, let $P =[p_{i,j}] = \big[ \Phi (e_{i,j}) \big] \in \M_N (\bofh)$
be the
positive Choi matrix determined by $\Phi$. Recall from the previous
lemma that the extended
$\Phi$ satisfies
$\Phi ( e_{i_1 i_1'} \otimes\dots\otimes e_{i_k i_k'}) = p_{i_1 i_1'} \cdots
p_{i_k i_k'}$. Thus 
\begin{align*}
\ip{x}{x}_\Phi &= \sum_{k,l \geq 0} 
\;\sum_{\substack{|w| = k\\|w'| =l}} \ip{w\otimes h_w}{w' \otimes
h_{w'}}_\Phi \\
&= \sum_{k\geq 0} 
\;\sum_{|w| = k} \ip{w\otimes h_w}{w' \otimes h_{w'}}_\Phi \\
&= \sum_{k\geq 0} 
\;\sum_{1\leq j \leq k} 
\;\sum_{1\leq i_j, i_j' \leq N} \ip{h_{i_1\cdots i_k}}{
\left(p_{i_1 i_1'} \cdots p_{i_k i_k'}\right) h_{i_1' \cdots i_k'}} .
\end{align*}
However, if we let $z_k = (h_w )_{|w| = k} \in \H^{(N^k)}$,
this quantity becomes
\[
\ip{x}{x}_\Phi = \sum_{k\geq 0} \ip{z_k}{P^{\otimes k} z_k}_{\H^{(N^k)}}
\geq 0.
\settowidth{\qedskip}{$\displaystyle 
\ip{x}{x}_\Phi = \sum_{k\geq 0} \ip{z_k}{P^{\otimes k} z_k}_{\H^{(N^k)}}
\geq 0.
$}
\addtolength{\qedskip}{-\textwidth}
\rlap{\hbox to-0.5\qedskip{\hfil\qed}}
\]
\renewcommand{\qed}{}
\end{proof}

\begin{defn}
Let $\N_\Phi = \{ x\in \T_N(\H) \mid \ip{x}{x}_\Phi = 0 \}$
be the kernel of the form
$\ip{\,\cdot\,}{\,\cdot\,}_\Phi$.
Since every positive semi-definite form satisfies the Cauchy-Schwarz
inequality,
$\N_\Phi$ is a subspace. We define the
\emph{Fock space of} $\Phi$ (\emph{or} $P = \left[ \Phi
(e_{i,j}) \right]$) \emph{over} $\H$ to be the Hilbert space completion
\[
\F_N (\H, \Phi) = \overline{ \T_N(\H) / \N_\Phi}^{\ip{\,\cdot\,}{\,\cdot\,}_\Phi}.
\]
\end{defn}

\begin{note}
The inner product on $\F_N (\H, \Phi)$ is given by 
\[
\ip{x+ \N_\Phi}{y + \N_\Phi}  =  
\ip{x}{y}_\Phi 
..
\]
We refer to the space
$e\otimes\H+\N_\Phi$
as the \emph{vacuum space} of $\F_N (\H, \Phi)$.
Recall
from the definition of the inner product
that
orthogonality is
preserved at
the level of the vacuum space,
unlike perhaps for words of larger length.
Hence there is no ambiguity in identifying it
with $\H$.
Further, notice
that in the
notation $\F_N (\H, \Phi)$,
reference to the
space $\H$ is really
redundant, since
it is fixed when
$\Phi$ is given. In
other words, the
construction is
totally determined
by the completely
positive map
\textup{(}equivalently, by the
associated positive
matrix\/\textup{)}.
\end{note}

The kernel $\N_\Phi$ can be explicitly identified in terms of $P$,
in fact the kernel of $P$. This
is implicit in the previous
proof, as is a Fourier expansion for vectors in
$ \F_N (\H, \Phi)$. This all follows from the
existence of projections onto `words of different lengths'. 
For $k\geq 0$, let $P_k$ be the map defined on finite sums
$x = \sum_{w\in\bbF_N^+} w 
\otimes h_w + \N_\Phi$
(in other words, finitely
many $h_{w}$ are nonzero)
by 
\[
P_k x = \sum_{|w| = k} w  \otimes h_w + \N_\Phi.
\]

\begin{lem} 
The maps $P_k $, for $k \geq 0$, extend to projections on
the Fock space $\F_N (\H, \Phi)$ with pairwise
orthogonal ranges. Further, we have
$
I =
\sideset{}{^{\smash{\oplus}}}{\textstyle\sum}\limits_{k\geq 0}
P_k$,
where the infinite sum is the limit in the strong operator topology.
\end{lem}

\begin{proof}
Each $P_k$ is clearly an idempotent.
Let $k \geq 0$ be fixed and let $x =  \sum_{w\in\bbF_N^+}
w  \otimes h_w + \N_\Phi$ be a finite sum. Then by orthogonality,
\begin{align*}
\left\| P_k x\right\|^2 &=  \sum_{|w| = k
= |w'|} \ip{w\otimes h_w}{w' \otimes h_{w'}} \\
     &\leq  \sum_{k \geq 0}  
\;\sum_{|w| = k = |w'|} \ip{w\otimes h_w}{w' \otimes h_{w'}} =
\left\|x\right\|^2. 
\end{align*}
Hence, $P_k$ extends to a contractive idempotent on $\F_N (\H, \Phi)$,
and as such, $P_k$ is a
(self-adjoint) projection on $\F_N (\H, \Phi)$.
Since the ranges of the $P_k$ are pairwise orthogonal, the
strong-operator-topology
limit $\sideset{}{^{\smash{\oplus}}}{\textstyle\sum}\limits_{k\geq 0}
P_k$
exists. However, this operator acts as the identity on a dense subset of
$\F_N (\H, \Phi)$, so it 
{\it is} in fact the identity operator.
\end{proof}

\begin{cor}
Every vector $x$ in $\F_N (\H,\Phi)$ has a representation of the form
\[
x = \sum_{k\geq0}P_{k}x=\sum_{w\in \bbF_N^+} w \otimes h_w + \N_\Phi.
\]
This representation is unique up to
choice of
the vectors
\[
P_k x = \sum_{|w| = k} w \otimes h_w + \N_\Phi.
\]
\end{cor}

More can be said in the Cuntz--Toeplitz case:
in unrestricted Fock space
vectors have \emph{bona fide} Fourier expansions,
since the vectors corresponding to words form an
orthonormal basis. In our  setting this can be seen as a relic of
$P= I_N$ having no kernel in that
case (This example is discussed further below.).
More generally, the positive matrix $P$ will have
nontrivial kernel, thus limiting the uniqueness of the Fourier expansion up to
representations
of the vectors $P_k x$.
We can obtain a tight upper bound on the norms of such vectors.

\begin{prop}\label{vector_estimate}
Let $h_w \in \H$ for each word
$w$ in $\bbF_N^+$ with
$|w| = k$. Then
\[
\left\| \sum_{|w|=k} w\otimes h_w + \N_\Phi \right\|^{2}
\leq \left\| P\right\|^{k} \left( \sum_{|w|=k} \left\|h_w\right\|^{2}
\right) .
\]
Further, this estimate is best possible in the sense that it can
always be asymptotically attained for some choice of
vectors $h_w$.
\end{prop}

\begin{proof}
Let $z = (h_w)_{|w|=k} \in \H^{(N^k)}$.
Then as in the proof of Theorem~\ref{possemidef} we
have
\begin{align*}
 \left\| \sum_{|w|=k} w\otimes h_w + \N_\Phi \right\|^2
&= \sum_{|w| = k = |w'|} \ip{w\otimes h_w}{w' \otimes h_{w'}} \\
&= \sum_{j=1}^k 
\;\sum_{1\leq i_j, i_j' \leq N} \ip{h_{i_1\cdots i_k}}
{\left(p_{i_1i_1'}\cdots 
p_{i_ki_k'}\right) h_{i'_1\cdots i'_k}} \\
&= \ip{z}{P^{\otimes k} z}   \leq   \left\|P\right\|^k \left\|z\right\|^2. 
\end{align*}
This establishes the desired inequality.
It is best possible since the vectors $z\in \H^{(N^k)}$ can
be chosen to approximate the norm of $P^{\otimes k}$. 
\end{proof}

These projections also allow us to illustrate further
the dependence of the Fock spatial structure
on the matrix $P$ in that they lead to a lucid identification of the kernel.

\begin{thm}
\label{ThmS:Fock.lucid}
The kernel $\N_\Phi$ is the closed span of the
pairwise orthogonal subspaces $P_k \N_\Phi$, for 
$k\geq 1$, given by
\begin{align*}
P_k \N_\Phi &= \left\{
\sum_{|w| = k} w\otimes h_w \Bigm| (h_w)_{|w| =k} \in\ker P^{\otimes k}
\right\} \\
&= \left\{ \sum_{|w| = k} w\otimes h_w \Bigm| (h_w)_{|w| =k} \simeq x_1 \otimes
\dots \otimes x_k \text{\textup{ with}} \right. \\  
&\qquad\qquad\qquad\qquad \left.
\vphantom{\sum_{|w| = k} w\otimes h_w \Bigm| (h_w)_{|w| =k} \simeq x_1 \otimes
\ldots \otimes x_k \text{\textup{ with}} }
x_1, \ldots, x_k  \in \H^{(N)}
\text{\textup{ and  some }} x_i \in \ker P \right\}.
\end{align*}
In particular, $\N_\Phi$ is completely determined by the kernel of $P$.
\end{thm}

\begin{proof}
With the existence of the projections $P_k$ proved in the previous
lemma, we have 
$\N_\Phi = 
\sideset{}{^{\smash{\oplus}}}{\textstyle\sum}\limits_{k\geq 0}
P_k \N_\Phi$. The characterization of the subspaces $P_k
\N_\Phi$ in terms of $P$ comes as a direct consequence of the proof of
Theorem~\ref{possemidef}. 
Indeed, from the end of that proof we see that $P_k \N_\Phi$
is determined by
$\ker P^{\otimes k}$, and this
subspace has the desired form.
\end{proof}

We finish this section by observing the ways in which
unrestricted Fock space is captured by our construction.

\begin{eg}
\label{ExaS:Fock.11}
Let
$P=I_N$ be the identity matrix in
$\M_N\left(\mathbb{C}\right)$. Then $P$ is the Choi matrix for the
completely positive map $\Phi\colon \M_N\rightarrow\mathbb{C}$ defined by
$\Phi\left(e_{i,j}\right)=1$ if $i=j$ and $0$ if $i\neq j$.
The corresponding Fock space $\fockspaceN\left(\mathbb{C},I_N\right)$ is the
standard unrestricted version. Indeed, since the kernel of $P$
is trivial, the null set $\N_\Phi=\left\{0\right\}$. Further,
for words $w=i_1\cdots i_k\in\bbF_N^+$ and letters $1\leq i,j\leq N$ we have
\[
\ip{iw\otimes\Omega}{jw\otimes\Omega}_\Phi
=\ip{\Omega}{\left(p_{i,j}p_{i_1,i_1}\cdots p_{i_k,i_k}\right)\Omega}.
\]
It follows that the vectors $\left\{w\otimes\Omega\right\}_{w\in\bbF_N^+}$ form
an orthonormal basis for the space, and structurally $\K =
\ell^2\left(\bbF_N^+\right)$
is obtained simply by identifying basis vectors $w\otimes \Omega$ with
$w$ for $w\in\bbF_N^+$ (see Example 4.1).

The construction also yields $\K = \ell^2\left(\bbF_N^+\right)$ through
what can be
considered as the biorthogonal setting. We will say more
about this at the end of the next section.
Let $E$ be the matrix in $\M_2\left(\mathbb{C}\right)$ with a one in each entry
and let $P=I_N\otimes E$.
Then the construction of $\fockspacetwoN\left(\mathbb{C},P\right)$ yields
$\K = \ell^2\left(\bbF_N^+\right)$
once again.
In this case, a nontrivial kernel for the matrix $P$ leads
to a nontrivial null set $\N_\Phi$. In particular, for words $w$ in
$\bbF_{2N}^+$
and letters $i\in\left\{1,\dots,2N\right\}$ with $1\leq i\leq N$, the vectors
\[
iw\otimes\Omega-\left(i+N\right)w\otimes\Omega
\]
belong to $\N_\Phi$, since $p_{i,i}=p_{i+N,i+N}=p_{i+N,i}=p_{i,i+N}=1$.
This collapse, together with the other entries of $P$,
shows that the vectors $\left\{w\otimes\Omega+\N_\Phi\right\}_{w\in\bbF_N^+}$
form an orthonormal basis for $\fockspacetwoN\left(\mathbb{C},P\right)$.
In other
words,  $\K = \ell^2\left(\bbF_N^+\right)$ is obtained structurally, with 
an obvious identification of orthonormal bases.
A similar analysis also shows that the construction
for $\fockspacetwoN\left(\H,I_{kN}\otimes E\right)$,
where $\H$ is $k$-dimensional
Hilbert space, yields $\ell^2\left(\bbF_N^+\right)^{\left(k\right)}$.
\end{eg}

\section{Creation operators}\label{S:creationops}







The Fock spaces from the previous section yield creation operators which
reduce to the
Cuntz--Toeplitz isometries
in the unrestricted
Fock space setting.

\begin{defn}
\label{DefS:creationops.1}The \emph{left creation operators}
$T=\left(T_{1},\dots,T_{N}\right)$ on
$\F_N (\H, \Phi)$ are linear transformations defined by
\[
T_{i}\left(w\otimes h+\N_\Phi\right)=\left(iw\right)\otimes h+\N_\Phi,
\]
where once again the product $iw$ is the free semigroup $\bbF_N^+$. 
\end{defn}

Since there is nontrivial null space in general,
we must check that these operators are well-defined.

\begin{prop}
\label{ProS:creationops.well}The operators $T=\left(T_1,\dots,T_N\right)$
are well-defined since $T_i\N_\Phi\subseteq\N_\Phi$ for $1\leq i\leq N$.
\end{prop}

\begin{proof}
It suffices to check that $T_iP_k\N_\Phi\subseteq P_{k+1}\N_\Phi$.
Let $x=P_kx=\sum_{\left|w\right|=k}w\otimes h_w+\N_\Phi$ belong to $\N_\Phi$,
and put $z=\left(h_w\right)_{\left|w\right|=k}\in
\H^{\left( N^k\right) }$.
Then $z\in\ker P^{\otimes\,k}$ by Theorem \ref{ThmS:Fock.lucid}, and we have
\[
\ip{T_ix}{T_ix}=\sum_{\left|w\right|=k=\left|w'\right|}
\ip{iw\otimes h_w}{iw'\otimes h_{w'}}_\Phi
=\ip{z}{p_{i,i}^{\left(N^k\right)}P^{\otimes\,k}z}=0.
\]
Thus $T_i$ belongs to $\N_\Phi$, and it follows that
$T_i$ is well-defined.
\end{proof}

For ease of presentation we shall suppress 
reference to the kernel $\N_\Phi$ for the rest of
this section.

\begin{rem}
\label{RemBefore4.6MarkII}
Our notion of the
creation operators is in principle similar to,
but yet quite different from, others
in the literature. One instance of this notion is the one
used in \cite{MuSo98} and \cite{MuSo99} in the
construction of covariant representations of
tensor algebras, such as
the Cuntz--Pimsner algebra. The
basic concepts in \cite{MuSo98} and \cite{MuSo99} are an inner
product $\ip{\,\cdot \,}{\,\cdot \,}$ on a module $E$
taking values in a $\ca$-algebra
$\fA$, and a completely positive mapping
$\Phi\colon\fA\rightarrow\BH$ where $\H$ is a Hilbert space. Then an extended
inner product is defined on a tensor algebra
over $E$, starting with $E\otimes\H$, as follows:
\[
\ip{a\otimes\xi}{b\otimes\eta}_{\operatorname*{new}}
:=\ip{\xi}{\Phi\left( \left\langle a\mid b\right\rangle\right) \eta},
\text{\qquad for }a,b\in E\text{ and }\xi,\eta\in\H.
\]
However, our construction and setting
are somewhat different. While we do start with
a completely positive map $\Phi\colon\M_{N}\rightarrow\BH$, our
construction yields a \emph{bona fide} inner product
which is defined by a certain natural extension
of the map $\tilde{\Phi}$ to the
`pre-$\uhfNI$' algebras $\M_{N^{k}}$
\textup{(}see Lemma \textup{\ref{LemTildePhi})}.
Further, in our most general case the creation
operators we get are not necessarily bounded.
Moreover, we need the
details from the intermediate steps of the construction.
\end{rem}

We also remark on the basis independence of this construction.

\begin{rem}
\label{RemBefore4.6}
The reader will notice that,
for a fixed completely positive map $\Phi$,
the notation $\fockspaceN\left(\H,\Phi\right)$
makes no reference to the basis for $\mathbb{C}^N$ used in our
construction. The reason for this is that the creation
operators obtained by using different orthonormal bases
are unitarily equivalent. To see this, let
$\left\{\xi_1,\dots,\xi_N\right\}$ and $\left\{\eta_1,\dots,\eta_N\right\}$
be orthonormal bases
for $\mathbb{C}^N$. Then the two Fock spaces constructed will be
spanned by vectors of the form, respectively,
$\xi_{i_1}\otimes\dots\otimes\xi_{i_k}\otimes h$ and
$\eta_{i_1}\otimes\dots\otimes\eta_{i_k}\otimes h$
(again, suppressing reference to the kernel). Let $U$ be the unitary
between these two spaces which identifies such spanning vectors \textup{(}it
is a unitary since the inner product is computed in the same way\/\textup{)}.
For $1\leq i\leq N$, let $T_i$ and $S_i$ be the $i$th creation operators
on the first and second space. Then for vectors
$x=\xi_{i_1}\otimes\dots\otimes\xi_{i_k}\otimes h$ we have
\begin{align*}
U^*S_iUx&=U^*S_i\left(\eta_{i_1}\otimes\dots\otimes\eta_{i_k}\otimes h\right)\\
&=U^*\left(\eta_i\otimes\eta_{i_1}\otimes\dots\otimes\eta_{i_k}\otimes  
h\right)
=T_ix.
\end{align*}
Thus, in this sense there is no ambiguity in using the
notation $\fockspaceN\left(\H,\Phi\right)$, given the completely positive
map $\Phi$.
\end{rem}

Through the behavior of the inner product, we can
describe the action of $T_{i}^{*}$ on the spanning vectors. Since
\[
\ip{T_{i}^{*}\left(e\otimes h\right)}{w\otimes h'}
=\ip{e\otimes h}{\left(iw\right)\otimes h'}=0,
\]
for all words $w\in\bbF_N^+$, we have $T_{i}^{*}\H=0$. Further
$\ip{T_{i}^{*}\left(w\otimes h\right)}{w'\otimes h'}=0$
unless $\left|w\right|=\left|w'\right|+1$. In this
case, if $w=i_{1}\cdots i_{k}$ and $w'=i'_{2}\cdots i'_{k}$, we have
\begin{align*}
\ip{T_{i}^{*}\left(w\otimes h\right)}{w'\otimes h'}
&=\ip{w\otimes h}{\left(iw'\right)\otimes h'}\\
&=\ip{h}{\left(p_{i_{1}i}p_{i_{2}i'_{2}}\cdots p_{i_{k}i'_{k}}\right)h'}.
\end{align*}

Let us summarize these actions in terms of our Fourier projections.

\begin{prop}
\label{ProS:creationops.2}For $1\leq i\leq N$, we have $P_{0}T_{i}=0$ and
$T_{i}P_{k}=P_{k+1}T_{i}$ for $k\geq0$.
\end{prop}

\begin{proof}
This simply follows from the above analysis. Indeed, it was observed that
$T_{i}^{*}$ annihilates $\ran P_{0}=\H$. Further, since
$I=
\smash[b]{
\sideset{}{^{\smash{\oplus}}}{\textstyle\sum}\limits_{j\geq 0}
}
P_{j}$ we have
\[
T_{i}P_{k}=P_{k+1}T_{i}P_{k}
=P_{k+1}T_{i}
\sideset{}{^{\smash{\oplus}}}{\sum}\limits_{j\geq 0}
P_{j}=P_{k+1}T_{i}.\settowidth
{\qedskip}{$\displaystyle T_{i}P_{k}=P_{k+1}T_{i}P_{k}
=P_{k+1}T_{i}
\sideset{}{^{\smash{\oplus}}}{\sum}\limits_{j\geq 0}
P_{j}=P_{k+1}T_{i}.$}\addtolength
{\qedskip}{-\textwidth}\rlap{\hbox to-0.5\qedskip{\hfil\qed}}
\]
\renewcommand{\qed}{}
\end{proof}

When the $T_{i}$ are isometries with pairwise orthogonal ranges,
the act of `pushing out' then
`pulling back' is given by $T_{i}^{*}T_{j}=\delta_{i,j}I$. In general
this action will not be so clean, since the
inner product twists each time $T_{i}$ or $T_{i}^{*}$ is applied.
Nonetheless, the action can be readily described in terms of the
matrix $P$.

\begin{thm}
\label{ThmS:creationops.3}Let $w\in\bbF_N^+$ and $h\in\H$.
Then for $1\leq i,j\leq N$ we have
\[
T_{i}^{*}T_{j}\left(w\otimes h\right)=w\otimes p_{i,j}h.
\]
\end{thm}

\begin{proof}
Since the vectors $w\otimes h$ span (not necessarily orthogonally)
\linebreak
the Fock space $\F_N (\H, \Phi)$, it suffices to examine
inner products
\linebreak
$\ip{T_{i}^{*}T_{j}\left(w\otimes h\right)}{w'\otimes h'}$
for $w,w'\in\bbF_N^+$ and $h,h'\in\H$.
If $w$ and $w'$ are words of different lengths, this inner product
is clearly zero. On the other hand, if $w=i_{1}\cdots i_{k}$ and
$w'=i'_{1}\cdots i'_{k}$ we have
\begin{align*}
\ip{T_{i}^{*}T_{j}\left(w\otimes h\right)}{w'\otimes h'}
&=\ip{\left(jw\right)\otimes h}{\left(iw'\right)\otimes h'}\\
&=\ip{h}{\left(p_{ji}p_{i_{1}i'_{1}}\cdots p_{i_{k}i'_{k}}\right)h'}\\
&=\ip{p_{i,j}h}{\left(p_{i_{1}i'_{1}}\cdots p_{i_{k}i'_{k}}\right)h'}\\
&=\ip{w\otimes p_{i,j}h}{w'\otimes h'}.
\end{align*}
The result now follows from the existence of the Fourier expansions for
vectors in $\F_N (\H, \Phi)$.
\end{proof}

Thus, at least formally, the transformations $T_{i}^{*}T_{j}$ can be
thought of as the tensor product of the identity on
unrestricted Fock space together with the operator $p_{i,j}$.
We can use this computation to discuss complete boundedness of $\Phi$ and
boundedness of the creation operators.

\begin{cor}
\label{CorS:creationops.4}When the representation of $\bbF_N^+$ on
$\F_N (\H, \Phi)$ determined by the creation operators
$T=\left(T_{1},\dots,T_{N}\right)$ yields
bounded operators, the map $\Phi$ is completely
bounded.
\end{cor}

\begin{proof}
Since $\Phi$ is completely positive and
$\Phi\left(e_{i,i}\right)=p_{i,i}$, we have
the inequality
\[
\left\|\Phi\right\|_{\operatorname*{cb}}
=\left\|\Phi\left(I_{N}\right)\right\|
=\left\|\sum_{i=1}^{N}\Phi\left(e_{i,i}\right)\right\|
\leq\sum_{i=1}^{N}\left\|p_{i,i}\right\|.
\]
However, from the theorem it follows that
$\left\|p_{i,i}\right\|=\left\|T_{i}^{*}T_{i}|_{\H}\right\|
\leq\left\|T_{i}\right\|^{2}<\infty$.
Whence, $\left\|\Phi\right\|_{\operatorname*{cb}}<\infty$.
\end{proof}

There is a partial converse of this result which
contains all the cases we are interested in.

\begin{cor}
\label{CorS:creationops.5}Let 
$T=\left(T_{1},\dots,T_{N}\right)$ be the
creation operators on
$\F_N (\H, \Phi)$, where $\Phi$ is completely
bounded with commutative range.
Then for $1\leq i\leq N$,
the norm of $T_{i}$ is given by
$\left\|T_{i}\right\|=\left\|p_{i,i}\right\|^{\frac{1}{2}}$.
\end{cor}

\begin{proof}
By Proposition \ref{ProS:creationops.2}, we have
$T_{i}^{*}T_{i}P_{k}=T_{i}^{*}P_{k+1}T_{i}
=P_{k}T_{i}^{*}T_{i}$ for $k\geq0$.
Hence $T_{i}^{*}T_{i}$ is diagonal with respect to the decomposition
$\F_N (\H, \Phi)=
\sideset{}{^{\smash{\oplus}}}{\textstyle\sum}\limits_{k\geq 0}
P_{k}\F_N (\H, \Phi)$, and as such
we obtain the norm identity
\[
\left\|T_{i}\right\|^{2}=\left\|T_{i}^{*}T_{i}\right\|
=\sup_{k\geq0}\left\|T_{i}^{*}T_{i}P_{k}\right\|.
\]
Let $k\geq0$ and let $h_{w}$ be vectors in $\H$ for each $\left|w\right|=k$.
If
$x=\sum_{\left|w\right|=k}w\otimes h_{w}$,
then from the theorem we have
\begin{align*}
\left\|T_{i}^{*}T_{i}x\right\|^{2}
&=\sum_{\left|w\right|=k=\left|w'\right|}
\ip{w\otimes p_{i,i}h_{w}}{w'\otimes p_{i,i}h_{w'}}\\
&=\ip{p_{i,i}^{\left(N^{k}\right)}z}{P^{\otimes k}p_{i,i}^{\left(N^{k}\right)}z},
\end{align*}
where $z=\left(h_{w}\right)_{\left|w\right|=k}$ belongs to
$\H^{\left(N^{k}\right)}$
and $p_{i,i}^{\left(N^{k}\right)}$ is the diagonal matrix in
$\M_{N^{k}}\left(\BH\right)$ with $p_{i,i}$ down the diagonal.
However, since the $p_{i,j}$ commute we arrive at the inequality
\[
p_{i,i}^{\left(N^{k}\right)}P^{\otimes k}p_{i,i}^{\left(N^{k}\right)}
\leq\left\|p_{i,i}\right\|^{2}P^{\otimes k}.
\]
This follows from the fact that if $R$ and $S$ are commuting
positive operators, then
\begin{equation}
0\leq RSR=\sqrt{S}R^{2}\sqrt{S}\leq\left\|R\right\|^{2}S.
\label{eqRSRBound}
\end{equation}
Thus we have
\[
\left\|T_{i}^{*}T_{i}x\right\|^{2}
\leq\left\|p_{i,i}\right\|^{2}
\ip{z}{P^{\otimes k}z}=
\left\|p_{i,i}\right\|^{2}
\left\|x\right\|^{2},
\]
so that
$\left\|T_{i}^{*}T_{i}P_{k}\right\|
\leq\left\|p_{i,i}\right\|$. Since
$T_{i}^{*}T_{i}P_{0}=p_{i,i}$,
it follows that 
$\left\|T_{i}^{*}T_{i}\right\|
=\left\|p_{i,i}\right\|$, as required.
\end{proof}

\begin{rem}
\label{RemCommutativityEssential}
The commutativity
assumption in Corollary \textup{\ref{CorS:creationops.5}} is essential.
The reader can check that, for the
case when $S$ in \textup{(\ref{eqRSRBound})} is a
projection, then the operator inequality
\textup{(\ref{eqRSRBound})} holds if \emph{and only if} $R$ and
$S$ are commuting.
\end{rem}

We
finish by pointing
out properties of the creation
operators which are naturally determined by
biorthogonal
wavelet representations.

\begin{cor}
\label{CorS:creationops.6}Let $S=\left(S_0,\dots,S_{N-1}\right)$
and $\tS =\left(\tS _0,\dots,\tS _{N-1}\right)$
form a biorthogonal
wavelet representation on $\H=L^2\left(\mathbb{T}\right)$
with invertible loop matrices $A$ and $\tA $.
Let $\S=\left[
\begin{array}{cc}
S & \tS 
\end{array}
\right]$
be a row matrix.
Let $P=\S^*\S$ be the
matrix in $\M_{2N}\left(\B\left(\H\right)\right)$ determined by the
representation,
as in Remark \textup{\ref{RemAfterLemS:wavelets.7}}. 
Let
\[
T=\left(T_1,\dots,T_{N},\tT_1,\dots,\tT_{N}\right)
\]
be the creation operators
acting on $\fockspacetwoN\left(\H,P\right)$. Then for $1\leq i,j\leq N$,
we have
\begin{enumerate}
\item $T_i^*T_j|_\H=S_{i-1}^*S_{j-1}=\left(AA^*\right)_{i,j}$,
\item $\tT_{i}^*\tT_{j}|_\H=\tS _{i-1}^*\tS _{j-1}
=\left(AA^*\right)^{-1}_{i,j}$,
\item $\tT_{i}^*T_j|_\H=T_i^*\tT_{j}|_\H=
\left\{
\renewcommand{\arraystretch}{1.25}
\begin{array}{l}
I\text{, if }i=j,\\
0\text{, if }i\neq j.
\end{array}
\right.
$
\end{enumerate}
\end{cor}

\begin{rem}
\label{RemS:creationops.7}
The corollary yields an esthetically pleasing relationship
between the biorthogonal wavelet representations and the creation operators
they determine. Indeed, the work of the last two
sections shows that the $S,\tS $ system completely
determines the Fock-space structure, and the
actions of the creation operators.
However, we are still trying to get a handle on what
this all means for these representations.
It was observed in Example \textup{\ref{ExaS:Fock.11}}
essentially how to obtain the
Cuntz--Toeplitz isometries as they sit inside the
biorthogonal class. But this is a little bit
misleading, for, if the previous corollary is applied with
the matrix $P$ defined by a  representation of $\O_N$, then
the Cuntz--Toeplitz isometries with \emph{infinite multiplicity}
are obtained. Perhaps our construction,
when applied to the biorthogonal wavelet representations,
somehow yields the appropriate creation operators for
the representations repeated with infinite multiplicity?
Or maybe if one wishes to study creation operators
associated with these representations, the infinite
multiplicity setting is a necessity?
In any event, there are a number of open problems
related to this class of creation operators.
\end{rem}




\begin{acknowledgements}
We are very grateful to Brian
Treadway for expert typesetting
and to Prof.\ Paul Muhly for enlightening discussions.
The second author
would like to thank the Department of Mathematics
at the University of Iowa for  kind hospitality during his visit.
We also thank the referees for helpful comments. 
\end{acknowledgements}


\begin{tabbing}
{\it E-mail address}:xx\= \kill
\noindent {\footnotesize\it Address}:
\>{\footnotesize\sc Department of Mathematics}\\
\>{\footnotesize\sc University of Iowa}\\
\>{\footnotesize\sc Iowa City, IA\quad 52242}\\
\>{\footnotesize\sc USA}\\
\\
{\footnotesize\it E-mail addresses}: 
\>{\footnotesize\sf jorgen@math.uiowa.edu} \\
\>{\footnotesize\sf dkribs@math.uiowa.edu}
\end{tabbing}

\end{document}